\documentclass[11pt,twoside]{amsart}
\usepackage{amscd}

\theoremstyle{plain}
\newtheorem{theorem}{Theorem}[subsection]

\newcommand{\rk}{\operatorname{rk}}
\newcommand{\res}{\operatorname{res}}

\newcommand{\diag}{\operatorname{diag}}
\newcommand{\Chow}{\operatorname{Chow}}
\newcommand{\Conv}{\operatorname{Conv}}

\newcommand{\End}{\operatorname{End}}
\newcommand{\GL}{\operatorname{GL}}
\newcommand{\SL}{\operatorname{SL}}

\newcommand{\Hilb}{\operatorname{Hilb}}

\newcommand{\Hom}{\operatorname{Hom}}

\newcommand{\bA}{{\mathbb A}}

\newcommand{\bF}{{\mathbb F}}
\newcommand{\bG}{{\mathbb G}}

\newcommand{\bM}{{\mathbb M}}
\newcommand{\bN}{{\mathbb N}}
\newcommand{\bP}{{\mathbb P}}

\newcommand{\bR}{{\mathbb R}}
\newcommand{\bT}{{\mathbb T}}

\newcommand{\bZ}{{\mathbb Z}}

\newcommand{\bmu}{\mathbb{\mu}}

\newcommand{\cA}{{\mathcal A}}
\newcommand{\cB}{{\mathcal B}}
\newcommand{\cC}{{\mathcal C}}
\newcommand{\cD}{{\mathcal D}}

\newcommand{\cF}{{\mathcal F}}

\newcommand{\cL}{{\mathcal L}}
\newcommand{\cM}{{\mathcal M}}

\newcommand{\cO}{{\mathcal O}}

\newcommand{\cQ}{{\mathcal Q}}
\newcommand{\cR}{{\mathcal R}}

\newcommand{\cX}{{\mathcal X}}
\newcommand{\cY}{{\mathcal Y}}

\newcommand{\hG}{\widehat G}

\newcommand{\hFun}{\widehat {\operatorname{Fun}}}

\newcommand{\hbM}{\widehat \bM}

\newcommand{\tB}{\tilde B}
\newcommand{\tC}{\tilde C}
\newcommand{\tG}{\tilde{G}}
\newcommand{\tR}{\tilde{R}}

\newcommand{\tT}{\tilde{T}}

\newcommand{\tX}{\tilde{X}}
\newcommand{\tY}{\tilde Y}
\newcommand{\tZ}{\tilde Z}

\newcommand{\tf}{\tilde f}
\newcommand{\tlambda}{\tilde\lambda}
\newcommand{\tmu}{\tilde\mu}

\newcommand{\tGamma}{\tilde\Gamma}
\newcommand{\tLambda}{\tilde\Lambda}
\newcommand{\tpi}{\tilde{\pi}}

\newcommand{\tcX}{\tilde{\mathcal X}}
\newcommand{\tcY}{\tilde{\mathcal Y}}

\newcommand{\wG}{\widetilde G}

\DeclareMathSymbol{\curvearrowright}{\mathrel}{AMSb}{"79}
\DeclareMathSymbol\rightsquigarrow {\mathrel}{AMSa}{"20}
\DeclareMathSymbol\square {\mathord}{AMSa}{"03}
\DeclareMathSymbol{\ltimes}         {\mathbin}{AMSb}{"6E}
\DeclareMathSymbol{\nmid}           {\mathrel}{AMSb}{"2D}
\DeclareMathSymbol{\twoheadrightarrow}  {\mathrel}{AMSa}{"10}

\newcommand{\ratmap}{- \kern -3pt \to}

\newcommand{\onto}{\twoheadrightarrow}

\newcommand{\Cone}{\operatorname{Cone}}

\newcommand{\main}{\operatorname{main}}

\newcommand{\id}{\operatorname{id}}

\newcommand{\Proj}{\operatorname{Proj}}

\newcommand{\Stab}{\operatorname{Stab}}
\newcommand{\Spec}{\operatorname{Spec}}

\newcommand{\rank}{\operatorname{rank}}

\newcommand{\im}{\operatorname{im}}

\newcommand{\Aut}{\operatorname{Aut}}

\newcommand{\Sym}{\operatorname{Sym}}

\theoremstyle{plain}

\newtheorem{lemma}[theorem]{Lemma}
\newtheorem{corollary}[theorem]{Corollary}

\newtheorem{proposition}[theorem]{Proposition}

\theoremstyle{definition}

\newtheorem{definition}[theorem]{Definition}

\newtheorem{example}[theorem]{Example}

\newtheorem{remark}[theorem]{Remark}   
\newtheorem{remarks}[theorem]{Remarks}

\newtheorem{acknowledgements}{Acknowledgments} 

\theoremstyle{remark}

\newcommand{\irr}{\rm irr}
\newcommand{\Gr}{\operatorname{Gr}}
\newcommand{\oGr}{\overline{\operatorname{Gr}}}
\newcommand{\ud}{\underline{d}}
\newcommand{\ui}{\underline{i}}

\newcommand{\oOmega}{\overline{\Omega}}
\newcommand{\wcA}{\widetilde{\cA}}

\author{Valery Alexeev and Michel Brion}
\address{Department of Mathematics\\
University of Georgia\\
Athens, GA 30602, USA}
\email{valery@math.uga.edu}
\address{Institut Fourier, B. P. 74\\
38402 Saint-Martin d'H\`eres Cedex, France}
\email{Michel.Brion@ujf-grenoble.fr}

\begin{document}

\begin{abstract}
We introduce a notion of stable spherical variety which includes the
spherical varieties under a reductive group $G$ and their flat equivariant
degenerations. Given any projective space $\bP$ where $G$ acts linearly, we
construct a moduli space for stable spherical varieties over $\bP$, that is,
pairs $(X,f)$, where $X$ is a stable spherical variety and $f : X \to \bP$
is a finite equivariant morphism. This space is projective, and its
irreducible components are rational. It generalizes the moduli space of
pairs $(X,D)$, where $X$ is a stable toric variety and $D$ is an effective
ample Cartier divisor on $X$ which contains no orbit. The equivariant
automorphism group of $\bP$ acts on our moduli space; the spherical
varieties over $\bP$ and their stable limits form only finitely many
orbits. A variant of this moduli space gives another view to the
compactifications of quotients of thin Schubert cells constructed by
Kapranov and Lafforgue. 
\end{abstract}

\bibliographystyle{amsalpha}

\title{Stable spherical varieties and their moduli}

\date{May 31, 2005}

\maketitle 

\tableofcontents

\setcounter{section}{-1}

\section{Introduction}

The starting point of this work is the construction and study of a
moduli space for stable toric pairs in \cite{Al}, and its
generalization to stable reductive pairs in \cite{AB-I,AB-II}. 
Both spaces, as well as the space of semiabelic pairs of \cite{Al},
parametrize certain pairs $(X,D)$, where $X$ is a projective variety
(possibly reducible, but not too singular), and $D$ is an effective
ample Cartier divisor on $X$. In addition, an algebraic group $G$ acts on
$X$ with finitely many orbits, none of them being contained in $D$; in the
two first cases, the invertible sheaf $\cO_X(D)$ is $G$-linearized. 


In both cases, the first step in the construction of the moduli space
is to classify the ``stable $G$-varieties'' $X$ and their orbit
structure, in terms of combinatorial invariants. The 
\emph{stable toric varieties} turn out to be unions of toric varieties
glued along torus-invariant subvarieties, and hence the well-known
classification of toric varieties may be applied. In the case of 
\emph{stable reductive varieties} (which includes the $G\times
G$-equivariant compactifications of a reductive group $G$ and their flat
equivariant degenerations), the classification may be reduced to that of
stable toric varieties with additional symmetries, by considering the
closures of appropriate torus orbits. 


Toric varieties and their non-commutative analogs, reductive
varieties, are examples of \emph{spherical varieties}. These may be
defined as the normal projective varieties where our reductive group
$G$ acts with finitely many orbits, so that this finiteness is
preserved under any equivariant modification. In the present work, we
introduce a notion of \emph{stable spherical variety} and we show the
existence of a projective moduli space of stable pairs in this setting.


Actually, we find it more convenient to consider pairs $(X,f)$, where $X$
is a stable spherical variety and $f:X \to \bP = \bP(V)$ is a finite
equivariant morphism to the projectivization of an arbitrary $G$-module
$V$. We then say that $X$ is a \emph{stable spherical variety over}
$\bP$. Clearly, for any pair $(X,D)$, the translates $gD$, $g\in G$, span a
base-point-free linear system; hence they define a pair $(X,f)$. It follows
that the moduli spaces of stable pairs $(X,D)$ are, in fact, special cases
of those of stable pairs $(X,f)$. For example, in the case when $G$ is a
torus, the space of pairs $(X,D)$ of \cite{Al} with linearized $\cO_X(D)$
corresponds to a multiplicity-free module $V$. (We note, however, that
\cite{Al} also deals with the case when the sheaf $\cO_X(D)$ is not
linearized.)


For the bulk of the paper we work with a connected reductive group $G$ over
an algebraically closed field $k$ of characteristic zero.  As we will see in
Section~\ref{sec:generalizations}, this assumption can be weakened in several
ways: $k$ need not be algebraically closed, and when $G$ is a split torus
one can work over an arbitrary base scheme, for example $\Spec\bZ$.
Section~\ref{sec:generalizations} also contains another generalization: to
stable varieties over a closed $G$-invariant subscheme $Z$ of $\bP(V)$. In
the case when $Z$ is a grassmanian, our moduli space is related 
to the compactifications of quotients of thin Schubert cells
constructed by Kapranov and Lafforgue, as we explain in
Example~\ref{ex:Lafforgue}.


Our construction goes along similar lines as in \cite{Al,AB-I,AB-II},
but an essential difference is that the combinatorial classification
of spherical varieties is unknown in general (although there are many
partial results, and a promising program; this will be discussed in
more detail below). So we systematically resort to qualitative 
arguments. 


A key ingredient is the finiteness of spherical varieties over a fixed
projective space $\bP$, up to equivariant automorphisms of $\bP$
(Theorem \ref{finiteness}). This result is deduced from the finiteness
of spherical orbit types in a fixed $G$-variety, proved in 
\cite[Section 3]{AB-III}; this follows, in turn, from a vanishing
theorem of Knop \cite{Kn-III} that implies the local rigidity of a class
of spherical varieties.


These finiteness results, proved by non-effective geometric arguments,
should rather be deduced from the classification (in progress) of
spherical varieties. The latter has been established by Luna and Vust
for varieties having a prescribed open $G$-orbit, see \cite{LV} and
also \cite{Kn-I}. Such spherical embeddings are described in terms of
combinatorial objects called colored fans, which also determine the
orbit structure; they generalize the fans of toric geometry. 


The next aim is to classify the spherical homogeneous spaces. This was
achieved by Wassermann in \cite{Wa} for spaces of rank at most $2$,
and then by Luna for groups $G$ of type $A$, see \cite{Lu}. The
combinatorial invariants introduced there make sense, in fact, for any
reductive group $G$, and they allowed Bravi and Pezzini to extend
Luna's result to groups of type $D$ \cite{BrPe}. But these
developments do not yet suffice to establish, in full generality, the
finiteness results that we need. 


Here our main combinatorial invariant is the \emph{moment polytope} 
of a polarized spherical variety $(X,L)$, a rational convex polytope
which governs the representations of $G$ in the spaces of sections of
powers of the ample line bundle $L$. In Corollary
\ref{finitepolytopes}, we obtain a fundamental boundedness property:


\noindent
\emph{For any rational convex polytope $Q$, there are only finitely many
isomorphism classes of polarized spherical varieties with moment
polytope $Q$, and given any bounded set $K$, there are only finitely
many moment polytopes contained in $K$.}  


In the toric case, the moment polytopes are precisely the integral
convex polytopes, and their faces correspond to toric subvarieties. 
But for a non-abelian group $G$, our understanding of moment polytopes
is very incomplete: their characterization is an open problem, and
their faces correspond only to certain Borel-invariant subvarieties,
see \cite{GS}.


Another open problem is the combinatorial classification of spherical
varieties over a fixed projective space $\bP$. In our moduli space,
they form finitely many orbits under the natural action of the
equivariant automorphism group of $\bP$. The orbit closures may be
described in terms of the Luna--Vust classification of spherical
embeddings; this will be developed elsewhere. In all examples that we
know of, there is a unique open orbit. Equivalently, all spherical
varieties over $\bP$ with prescribed weight group and moment polytope
are degenerations of a unique one. 


This is closely related to the Delzant conjecture in symplectic
geometry (see \cite{De} and also \cite{Wo}), which asserts that any
compact multiplicity-free Hamiltonian manifold is uniquely determined
by its moment polytope and principal isotropy subgroup. It suggests
that the \emph{nonsingular} spherical varieties over $\bP$ are classified
by their weight group and moment polytope.


In this work, we have endeavoured to make the exposition
self-contained, in particular, independent of the classification of
spherical embeddings. Thus, we have provided complete proofs of some
statements (Propositions \ref{CM}, \ref{isotrivial} and \ref{model})
which generalize the corresponding results in \cite{AB-I,AB-II}, but
where the original proofs relied on the classification. Examples at
the end of Sections 1, 2, and 4 discuss stable toric and reductive
varieties in detail and exhibit some new features of stable spherical  
varieties, which form their natural and definitive generalization.

\begin{acknowledgements}
The first author would like to thank the \'Ecole Normale Sup\'erieure
for its hospitality; his research was also partially supported by NSF
under DMS-0401795. 
\end{acknowledgements}

\section{Varieties with reductive group action}

\subsection{Preliminaries}

In this section, we fix the notation and gather preliminary results
concerning algebraic groups and varieties. We use \cite{Ha} as a
general reference for algebraic geometry, and \cite{PV}, \cite{Gr} for 
algebraic transformation groups.

The ground field $k$ is algebraically closed, of characteristic zero. By a
\emph{variety}, we mean a connected separated reduced scheme of finite type
over $k$; by a \emph{subvariety}, we mean a closed subvariety. 

A variety $X$ equipped with an algebraic action of a linear algebraic
group $G$ is called a $G$-\emph{variety};  an equivariant morphism 
between $G$-varieties is called a $G$-\emph{morphism}.
 
We will only consider $G$-\emph{quasiprojective} varieties, i.e., 
those admitting an ample $G$-linearized invertible sheaf. 
A quasiprojective $G$-variety $X$ is $G$-quasiprojective whenever $X$
is normal (by \cite{Su}), or when $G$ is finite. 

Throughout this paper, we denote by $G$ a connected reductive  
algebraic group ; we choose a Borel subgroup $B$ of $G$, and a maximal
torus $T$ of $B$. Then $B=TU$, where $U$ denotes the unipotent part of
$B$. We denote by $W=N_G(T)/T$ the Weyl group of $(G,T)$, and by 
$B^- = T U^-$ the unique Borel subgroup of $G$ such that 
$B^- \cap B = T$. The \emph{weight group} of $G$ is the character
group $\cX(T)$, identified with the character group of $B$ and denoted
by $\Lambda$.

By a $G$-\emph{module}, we mean a rational, possibly infinite-dimensional
$G$-module. Since $G$ is reductive, any $G$-module is semi-simple.
Moreover, any simple $G$-module $V$ is finite-dimensional and contains 
a unique line of $B$-eigenvectors; the corresponding weight 
$\lambda \in \Lambda$ determines $V$ uniquely. This yields a bijection
$\lambda\mapsto V(\lambda)$ from the set $\Lambda^+$ of dominant weights,
to the set of isomorphism classes of simple $G$-modules.

We may view $\Lambda$ as a lattice in the real vector space 
$\Lambda_{\bR} := \Lambda \otimes_{\bZ} {\bR}$; then $\Lambda^+$ is
the intersection of $\Lambda$ with the positive Weyl chamber
$\Lambda^+_{\bR} := \bR_{\ge 0} \Lambda^+$. 

Given a $G$-module $V$, the \emph{weight set} $\Lambda^+(V)$ is the
set of dominant weights of simple submodules of $V$; the 
\emph{weight cone} $C(V) := \bR_{\ge 0} \, \Lambda^+(V)$ 
is the subcone of $\Lambda^+_{\bR}$ generated by $\Lambda^+(V)$. 
We say that $V$ is \emph{multiplicity-free} if it is the direct sum of
pairwise non-isomorphic simple $G$-modules. Then $V$ is uniquely
determined by its weight set.

Let $X$ be an affine $G$-variety. Then the \emph{affine ring}
$R := H^0(X,\cO_X)$ is a $G$-module; this defines the weight set 
$\Lambda^+(X) := \Lambda^+(R)$ and the weight cone 
$C(X) := C(R)$. Both are discrete invariants of the $G$-variety $X$,  
that satisfy the following basic properties (see \cite{Sj} for the
proofs).

\begin{lemma}\label{weightset}
(i) If $X$ is irreducible, then $\Lambda^+(X)$ is a finitely generated
submonoid of $\Lambda^+$. Thus, $C(X)$ is a rational polyhedral convex
cone.

\smallskip

\noindent
(ii) For an arbitrary $X$ and a $G$-subvariety $Y$ holds
$\Lambda^+(Y) \subseteq \Lambda^+(X)$, whence $C(Y)\subseteq C(X)$.
Moreover, $\Lambda^+(X)$ is the union of the weight monoids
$\Lambda^+(Y)$, where $Y$ runs over the irreducible components of $X$.

\smallskip

\noindent
(iii) Let $X'$ be an affine $G$-variety and $f:X' \to X$ a finite
surjective $G$-morphism. Then 
$\Lambda^+(X) \subseteq \Lambda^+(X') \subseteq 
\frac{1}{N}\Lambda^+(X)$ for some positive integer $N$. Thus, 
$C(X') = C(X)$.
\end{lemma}

The affine $G$-variety $X$ is \emph{multiplicity-free}, if its affine
ring is a multiplicity-free $G$-module. For an irreducible
$G$-variety, the multiplicity-freeness is equivalent to the existence
of a dense $B$-orbit, and also to the finiteness of the number of
$B$-orbits. In particular, every affine multiplicity-free $G$-variety
$X$ contains only finitely many $G$-orbits. Also, recall the following
result (see, e.g., \cite[Lemma 1.3]{AB-III} for the easy proof):

\begin{lemma}\label{diag}
The $G$-automorphism group of any affine multiplicity-free
$G$-variety is a diagonalizable linear algebraic group.
\end{lemma}

Next let $X$ be an irreducible (possibly non-affine) $G$-variety. Then
$G$, and hence $B$, acts on the function field $k(X)$. The set of
weights of $B$-eigenvectors in $k(X)$ is a subgroup of $\Lambda$
that we denote by $\Lambda(X)$ and call the 
\emph{weight group of $X$}; this is a birational invariant of the
$G$-variety $X$. The abelian group $\Lambda(X)$ is free of finite
rank: the \emph{rank of} $X$, denoted by $\rk(X)$. 

If $X$ is affine, then 
$\Lambda^+(X)\subseteq \Lambda(X) \cap C(X)$,
and $\Lambda^+(X)$ generates the group $\Lambda(X)$. Thus, the rank of
$X$ is the dimension of the cone $C(X)$. 

For an arbitrary irreducible $G$-variety $X$, we record the
following result, a consequence of Lemma \ref{weightset} 
together with the local structure of $G$-varieties 
\cite[Section 2.2]{Kn-II}.

\begin{lemma}\label{weightgroup}
(i) Let $Y$ be an irreducible $G$-subvariety of $X$. Then $\Lambda(Y)$
is a subgroup of $\Lambda(X)$. In particular, $\rk(Y) \le \rk(X)$.

\noindent
(ii) Let $X'$ be an irreducible $G$-variety and $f : X' \to X$ a finite
surjective $G$-morphism. Then $\Lambda(X')$ contains $\Lambda(X)$ as a
subgroup of finite index. In particular, $\rk(X') = \rk(X)$.
\end{lemma}

\subsection{Polarized varieties}

We now adapt classical notions related to polarized varieties (see,
e.g., \cite{Vi}) to our equivariant setting. 

\begin{definition}
A \emph{polarized $G$-variety} is a pair $(X,L)$, where $X$ is a
projective $G$-variety and $L$ is an ample $G$-linearized invertible
sheaf on $X$. 

A \emph{$G$-morphism} $\varphi:(X',L') \to (X,L)$ of polarized
$G$-varieties is a pair $(f: X' \to X, ~\gamma: L' \to f^*L)$,
where $f$ is a $G$-morphism and $\gamma$ is an \emph{isomorphism} 
of $G$-linearized sheaves.
\end{definition}

For any such pair $(f,\gamma)$, the morphism $f$ is finite,
since $f^* L$ is ample. Also, note that $\gamma$ is uniquely
determined by $f$ up to scalar multiplication, since
$H^0(X',\cO_{X'})=k$ (as $X'$ is reduced and connected). 

Thus, the $G$-automorphisms of $(X,L)$ are exactly the pairs 
$(f,z f^*)$, where $f$ is a $G$-automorphism of $X$ that fixes 
the isomorphism class $[L]$ of the $G$-linearized invertible sheaf
$L$, and $z$ is a non-zero scalar. This yields an isomorphism 
$$
\Aut^G(X,L) \simeq \bG_m \times 
\{ f \in \Aut^G(X) ~\vert~ f^*[L] = [L] \}.
$$
Moreover, the $G$-isomorphism classes of polarized varieties are the
equivalence classes of pairs $(X,[L])$ under isomorphisms of
$G$-varieties.

To any polarized $G$-variety $(X,L)$ is associated its 
\emph{section ring}
$$
R(X,L) := \bigoplus_{n=0}^{\infty} H^0(X,L^n),
$$
where $L^n$ denotes the $n$-th tensor power of $L$. This is a
finitely generated graded algebra on which $G$ acts by automorphisms,
via its natural action on every $H^0(X,L^n)$. In other words, $R(X,L)$
carries an action of the group
$$
\tG := \bG_m \times G,
$$ 
where the multiplicative group $\bG_m$ acts with weight $n$ on 
$H^0(X,L^n)$. Note that $\tG$ is a connected reductive group
with Borel subgroup 
$\tB := \bG_m \times B$, maximal torus $\tT := \bG_m \times T$, 
weight group $\tilde\Lambda := \bZ \times \Lambda$, and set of
dominant weights $\tLambda^+ := \bZ \times \Lambda^+$. Moreover, 
$$
\tX := \Spec \, R(X,L)
$$ 
is an affine $\tG$-variety, the \emph{affine cone} over $X$. This
variety has a special closed point $0$, associated with the maximal
homogeneous ideal of $R(X,L)$ and fixed by $\tG$. As is well-known
(and follows from \cite[Proposition I.2.3]{Ra}), the natural
map
$$
\pi:\tX \setminus \{0\}\to X = \Proj \, R(X,L)
$$  
is a principal $\bG_m$-bundle, and there are isomorphisms
of $G$-linearized sheaves $\cO_X(n)\simeq L^n$ for all $n \in \bZ$.

Also, note that $X$ is normal (resp. semi-normal) if and only if $\tX$
is (see \cite[Lemma 2.1]{AB-II} for the semi-normality). 

Given another polarized $G$-variety $(X',L')$ with affine cone $\tX'$,
one easily sees that the $G$-morphisms $(X',L') \to (X,L)$
are in bijective correspondence with the \emph{finite} 
$\tG$-morphisms $\tX'\to \tX$. In particular,
\begin{equation}\label{aut}
\Aut^G(X,L) \simeq \Aut^{\tG}(\tX).
\end{equation}
On the other hand, any $G$-subvariety $Y$ of $X$ yields a polarized
$G$-variety $(Y,L) := (Y,L\vert_Y)$, and hence a finite morphism
$\tY \to \tX$ which is injective, but need not be a closed
immersion. 
 
Next we define the \emph{weight set} 
$\tLambda^+(X,L) := \tLambda^+(\tX)$, the \emph{weight cone}
$\tC(X,L) := C(\tX)$, and the \emph{moment set} 
$$
Q(X,L) := \tC(X,L) \cap (\{1\}\times \Lambda^+_{\bR})
\subset \tLambda_{\bR} = \bR \times \Lambda_{\bR}.
$$
Then $\tC(X,L) = \Cone Q(X,L)$ (the cone over the moment set).
Moreover, $Q(X,L)$ is a finite union of rational convex polytopes 
associated with the irreducible components of $X$ (as follows 
from Lemma \ref{weightset}, see \cite{Sj} for details). 
We identify $Q(X,L)$ with its projection to $\Lambda^+_{\bR}$. 

If $X$ is irreducible, then $Q(X,L)$ is a rational convex polytope 
in $\Lambda^+_{\bR}$, the \emph{moment polytope}. 
Moreover, the first projection $\tLambda \to \bZ$
yields an exact sequence
$$
0 \to \Lambda(X) \to \tLambda(X,L) \to \bZ \to 0,
$$
where $\tLambda(X,L) :=\tLambda(\tX)$ is the \emph{weight group}
of $(X,L)$. The real vector space $\Lambda(X)_{\bR}$ is spanned by 
the differences of the points of $Q(X,L)$ (for these results, see 
\cite{Sj}). In particular, $\dim Q(X,L) = \rk(X)$.

Clearly, the weight set, weight cone, and moment set of the pair
$(X,L)$ depend only on its $G$-isomorphism class. Also, note the 
following consequence of Lemmas \ref{weightset}(iii) and
\ref{weightgroup}(ii):

\begin{lemma}\label{moment}
Let $(f,\gamma):(X',L')\to (X,L)$ be a morphism of polarized
$G$-varieties, where $f$ is surjective. Then $Q(X',L') = Q(X,L)$. 
\end{lemma}

We say that a polarized $G$-variety $(X,L)$ is \emph{multiplicity-free},
if every $G$-module $H^0(X,L^n)$ is multiplicity-free (as is easy to
see, it suffices to check the multiplicity-freeness of $H^0(X,L^n)$
for $n\gg 0$). Equivalently, the affine cone $\tX$ is a
multiplicity-free $\tG$-variety. Then (\ref{aut}) and Lemma \ref{diag} 
yield:

\begin{lemma}\label{diago}
The $G$-automorphism group of any multiplicity-free polarized
$G$-variety is a diagonalizable linear algebraic group.
\end{lemma}

Returning to general polarized varieties, the simplest examples are,
of course, the pairs $(\bP(V),\cO(1))$, where $V$ is a
finite-dimensional $G$-module with projectivization 
$\bP(V) := \Proj \Sym(V^*),$  and $\cO(1)$ is equipped with its
natural $G$-linearization. Then $R(\bP(V),\cO(1)) = \Sym(V^*)$. 
We now introduce a stronger notion of polarization by prescribing a
$G$-morphism to some $(\bP(V),\cO(1))$.

\begin{definition}
Let $V$ be a finite-dimensional $G$-module.
A $G$-\emph{variety over} $\bP(V)$ is a pair $(X,f)$, where $X$ is a
projective $G$-variety and $f: X \to \bP(V)$ is a \emph{finite}
$G$-morphism.

A \emph{$G$-morphism} from another pair $(X',f')$ to $(X,f)$ is a
$G$-morphism $\varphi: X' \to X$ such that $f' \circ \varphi = f$.
\end{definition}

The $G$-varieties $(X,f)$ over $\bP(V)$ are in bijective 
correspondence with the polarized $G$-varieties $(X,L)$ 
equipped with a $G$-module map 
$$
\gamma: V^* = H^0(\bP(V), \cO(1)) \to H^0(X,L)
$$ 
such that the image of $\gamma$ is base-point-free. Namely,
one associates to $(X,f)$ the sheaf $L := f^*\cO(1)$
and the map $\gamma := f^*$. This also yields a finite $\tG$-morphism
\begin{equation}\label{lift}
\tf: \tX \to \bA(V),
\end{equation}
where $\tX$ is the affine cone over $(X,L)$, and $\bA(V)$ denotes 
the affine space $\Spec \Sym(V^*)$. The group $\tG = \bG_m \times G$
acts on $\bA(V)$ via the scalar action of $\bG_m$ and the given action
of $G$.

Clearly, there are only finitely many $G$-morphisms 
$\varphi: (X',f') \to (X,f)$ between any two prescribed $G$-varieties
over $\bP(V)$. Every such morphism defines a (special) $G$-morphism
$(\varphi,\id): (X',L') \to (X,L)$ of polarized varieties. 

The $G$-automorphism group $\Aut^G \bP(V)$ acts on the set of 
$G$-varieties $(X,f)$ over $\bP(V)$, and the isomorphism class $[L]$
of the $G$-linearized sheaf $L=f^*\cO(1)$ depends only on the
orbit. In the multiplicity-free case, this observation may be refined
as follows:

\begin{lemma}\label{class}
Given a finite-dimensional $G$-module $V$, there is a bijection
between:

\smallskip

\noindent
the $\Aut^G \bP(V)$-orbits of multiplicity-free $G$-varieties $(X,f)$
over $\bP(V)$, and

\smallskip

\noindent
the triples $(X,[L],F)$, where $(X,L)$ is multiplicity-free and $F$
is a subset of the weight set of $V^*$ such that $L$ is generated by a
submodule of global sections with weight set $F$. 
\end{lemma}

\begin{proof}
Consider a pair $(X,f)$ and the associated map 
$f^* : V^* \to H^0(X,L)$. Since the $G$-module $H^0(X,L)$ is
multiplicity-free, the image of $f^*$ is uniquely determined by its
weight set $F$. Moreover, the triple $(X,[L],F)$ only depends on the
$\Aut^G \bP(V)$-orbit of $(X,f)$.

Conversely, given a triple $(X,[L],F)$, let $H^0(X,L)_F$ be the unique 
$G$-submodule of $H^0(X,L)$ with weight set $F$. Our assumptions on
weight sets imply the existence of a surjective $G$-module map 
$\gamma : V^* \to H^0(X,L)_F$. Moreover, any two such maps are conjugate
under an element of $\GL(V^*)^G \simeq \GL(V)^G$. Finally, any
isomorphism of $G$-linearized sheaves $L \to L'$ yields a $G$-module
isomorphism $H^0(X,L') \to H^0(X,L)$ which maps isomorphically
$H^0(X,L')_F$ to $H^0(X,L)_F$. Thus, the corresponding maps
$\gamma$, $\gamma'$ are still conjugate under $\GL(V^*)^G$.
\end{proof}

The structure of $\Aut^G \bP(V)$ is easily described: since $G$ is
connected, one obtains an exact sequence 
$$
1 \to \bG_m \to \GL(V)^G \to \Aut^G \bP(V) \to 1,
$$
where $\GL(V)^G$ denotes the group of linear $G$-automorphisms of
$\bA(V)$, and $\bG_m$ acts on $\bA(V)$ by scalar
multiplication. Further, we have a canonical isomorphism of
$G$-modules
\begin{equation}\label{dec}
V \simeq \bigoplus_{\lambda \in F} E(\lambda) \otimes V(\lambda),
\end{equation}
where $F$ is a finite subset of $\Lambda^+$, and each 
$E(\lambda) = \Hom^G(V(\lambda),V)$ is a non-zero vector space of
finite dimension. Here, of course, $G$ acts on every summand
$E(\lambda) \otimes V(\lambda)$ via its action on $V(\lambda)$. Then 
\begin{equation}\label{prod}
\GL(V)^G \simeq \prod_{\lambda \in F} \GL(E(\lambda)),
\end{equation}
where each factor $\GL(E(\lambda))$ acts on $V$ via its natural
action on $E(\lambda)$. 

\subsection{Spherical varieties}

Recall that an affine irreducible $G$-variety $X$ is \emph{spherical}
if it is normal and multiplicity-free. Likewise, a polarized variety
$(X,L)$ is spherical if $X$ is normal and $(X,L)$ is
multiplicity-free; equivalently, the affine cone $\tX$ is a spherical
$\tG$-variety. We then say for brevity that $(X,L)$ is a PSV.

We now gather some fundamental facts on polarized spherical varieties.

\begin{proposition}\label{psv}
Let $(X,L)$ be a PSV with weight group $\Gamma$ and moment polytope
$Q$, and let $Y$ be an irreducible $G$-subvariety. Then:

\smallskip

\noindent
(i) $(Y,L)$ is a PSV.

\smallskip

\noindent
(ii) $\tilde\Lambda(Y,L)$ is a direct summand of $\Gamma$. Thus,
$\Lambda(Y)$ is a direct summand of $\Lambda(X)$.

\smallskip

\noindent
(iii) $Q(Y,L)$ is a face of $Q(X,L)$, which determines $Y$ uniquely. 

\smallskip

\noindent
(iv) The restriction map $H^0(X,L^n) \to H^0(Y,L^n)$ is surjective for
all $n\ge 0$. 

\smallskip

\noindent
(v) As $G$-modules,
$$
H^0(X,L^n) \simeq \bigoplus_{\lambda\in\Lambda^+,\;
(n,\lambda)\in \Gamma \cap nQ} V(\lambda).
$$
Equivalently, as $\tG$-modules,
$$
R(X,L) \simeq \bigoplus_{\tlambda \in \Gamma \cap \Cone(Q)}
V(\tlambda).
$$

\smallskip

\noindent
(vi) $L$ is globally generated.

\smallskip

\noindent
(vii) $H^i(X,L^n)=0$ for all $i\ge 1$, $n\ge 0$.
\end{proposition}

\begin{proof}
The natural map $f:\tY \to \tX$ between affine cones is
finite and injective. Thus, $f(\tY)$ is an irreducible
$\tG$-subvariety of the affine spherical variety $\tX$. Hence $f(\tY)$
is normal (e.g., by \cite[Lemma 2.2]{AB-I}). It follows that 
$f$ is an isomorphism; this implies (i) and (iv). 

Likewise, (ii), (iii) and (v) follow from their affine analogs proved
for example in [loc. cit.].

For (vi), if $L$ is not globally generated, then its base locus
contains a closed $G$-orbit $Y$. Since $L$ is ample and $Y$ is a flag
variety, then $H^0(Y,L)\ne 0$ which contradicts (iv).

Finally, (vii) is well-known, see, e.g., \cite[Corollary 5.8]{AB-II}. 
\end{proof}

\begin{lemma}\label{cover}
Let $X$ be a spherical $G$-variety, $Y$ an irreducible $G$-variety,
and $f: X \to Y$ a finite surjective $G$-morphism. Then:

\smallskip

\noindent
(i) $f^{-1}(Z)$ is irreducible for every irreducible $G$-subvariety
$Z$ of $Y$. Equivalently, the preimage of every $G$-orbit is
a unique $G$-orbit. 

\smallskip

\noindent
(ii) The pair $(X,f)$ is uniquely determined (up to $G$-isomorphism
over $Y$) by the data of $Y$ and the weight group $\Lambda(X)$. In
particular, if $\Lambda(X) = \Lambda(Y)$, then $f$ is the
normalization map.
\end{lemma}

\begin{proof}
(i) is a direct consequence of Proposition \ref{psv} (iii) combined
with Lemma \ref{moment}.

(ii) First we consider the case where $Y = G/I$ is a unique
$G$-orbit. Then we may write $X = G/H$, where $H$ is a subgroup of
finite index of $I$. Thus, the connected component $I^0$ satisfies
$I^0 \subseteq  H \subseteq I$. Now, by \cite[Section 5.2]{BP}, we
may choose $B$ and $T$ so that $BI^0$ is open in $G$, and 
$I = I^0 (T\cap I)$. Thus, $H = I^0 (T\cap H)$. Moreover, 
$\Lambda(X) = \cX(T/T\cap H)$ as a subgroup of $\Lambda=\cX(T)$, by 
\cite[Section 2.9]{BP}. So $\Lambda(X)$ determines uniquely the
subgroup $T\cap H$ of $T$ which, in turn, determines uniquely $H$.

In the general case, let $G/I$ be the open orbit in $Y$. Since $f$ is
finite, $f^{-1}(G/I) = G/H$, where $H$ is as above. Moreover, $X$ is
the normalization of $Y$ in the function field $k(G/H)$.
\end{proof}

Our final preliminary result concerns the $G$-automorphism group of a
spherical variety $X$. This group acts on the function field $k(X)$
and preserves each subset $k(X)^{(B)}_{\lambda}$ of 
$B$-eigenvectors with prescribed weight $\lambda\in\Lambda(X)$. Since
$X$ contains a dense $B$-orbit, such an eigenvector is determined by 
its weight up to multiplication by a non-zero scalar. This defines 
a character $\chi_{\lambda}: \Aut^G(X)\to \bG_m$ and, in turn, 
a homomorphism
$$
\iota = \iota_X : \Aut^G(X) \to \Hom(\Lambda(X), \bG_m),
\quad \varphi \mapsto (\lambda \mapsto \chi_{\lambda}(\varphi)),
$$
where the target is the torus with character group $\Lambda(X)$.
Now \cite[Theorems 5.1, 5.5]{Kn-IV} yields the following

\begin{lemma}\label{automorphisms}
Let $X$ be a spherical variety. Then every irreducible
$G$-subvariety $Y$ is invariant under $\Aut^G(X)$, and the diagram
$$
\CD
\Aut^G(X) @>{\iota_X}>>  \Hom(\Lambda(X),\bG_m) \\
@V{\res_Y}VV  @V{\res_{\Lambda(Y)}}VV \\
\Aut^G(Y) @>{\iota_Y}>>  \Hom(\Lambda(Y),\bG_m) \\
\endCD
$$
commutes. Moreover, $\iota_X$ and $\iota_Y$ are injective.
\end{lemma}

This lemma is easily checked directly if $X$ is affine, see, e.g. 
\cite[Lemma 1.3]{AB-III}. By (\ref{aut}), it also applies to
$\Aut^G(X,L)$, where $(X,L)$ is any polarized spherical variety.

We now come to a key finiteness result.

\begin{theorem}\label{finiteness}
Let $V$ be a finite-dimensional $G$-module. Then there are only
finitely many orbits of spherical varieties over $\bP(V)$, for the
action of $\Aut^G \, \bP(V)$.
\end{theorem}

\begin{proof}
Consider such a variety $(X,f)$ and the associated finite
$\tG$-morphism $\tf : \tX \to V$ of (\ref{lift}).
By \cite[Corollary 3.3]{AB-III}, there are only finitely many
possibilities for the image $\tf(\tX)$, up to the action of 
$\GL(V)^G$. So we may fix $\tY := \tf(\tX)$. 

By Lemma \ref{weightgroup}, $\tLambda(\tX)$ contains 
$\tilde\Lambda(\tY)$ as a subgroup of finite index. 
Since the torsion subgroup of the quotient 
$\tLambda/\tLambda(\tY)$ is finite,
there are only finitely many possibilities for the weight group 
$\tilde\Lambda(\tX)$. Moreover, by Lemma \ref{cover}, $\tY$ and the
latter group determine uniquely the pair $(\tX,\tf)$, and hence $(X,f)$.
\end{proof}

Together with Lemma \ref{class}, this implies readily the following:

\begin{corollary}\label{finitequotients}
Let $V$ be a finite-dimensional $G$-module. Then there are only
finitely many $G$-isomorphism classes of PSV's $(X,L)$, where $L$ is
generated by a $G$-module of global sections which is a quotient of
$V$. 
\end{corollary}

This implies, in turn, another finiteness result:

\begin{corollary}\label{finitepolytopes}
Let $K$ be a bounded subset of $\Lambda_{\bR}$. Then there 
are only finitely many $G$-isomorphism classes of PSV's with moment
polytope contained in $K$. 
\end{corollary}

\begin{proof}
Let $(X,L)$ be a PSV such that $Q(X,L) \subseteq K$. By Proposition
\ref{psv}, $L$ is generated by a $G$-module of global sections which
is a quotient of the $G$-module 
$\bigoplus_{\lambda\in \Lambda^+ \cap K} V(\lambda)^*$.  
So the statement follows from Corollary \ref{finitequotients}.
\end{proof}

In particular, any prescribed bounded set contains only finitely many 
moment polytopes of PSV's. This boundedness property is not obvious,
since the polytopes under consideration may have non-integral  
vertices; see, e.g., Examples 1.4(2), 1.4(3). 

\subsection{Examples} 

We now consider three natural subclasses of the class of polarized
spherical varieties, beginning with the simplest one:

\medskip

\noindent
1) \emph{Toric varieties.} 

Here $G=T$ is a torus with character group $\Lambda$. The spherical
varieties for $T$ are just the normal varieties where $T$ acts with
a dense orbit, i.e., the toric varieties for a quotient of $T$.  

As is well-known, the assignment $X \mapsto (\Lambda(X),C(X))$ yields
a bijective correspondence from the affine toric varieties for a
quotient of $T$, to the pairs $(\Gamma, C)$, where $\Gamma$ is a
subgroup of $\Lambda$, and $C$ is a rational polyhedral convex cone in
$\Gamma_{\bR}$ with non-empty interior. Moreover, 
$\Aut^T(X) = \Hom(\Gamma,\bG_m)$. 

Likewise, the assignment $(X,L)\mapsto (\tLambda(X,L),Q(X,L))$ yields
a bijective correspondence from the polarized toric varieties for a 
quotient of $T$, to the pairs $(\Gamma,Q)$, where $Q$ is a convex
lattice polytope in $\{1\}\times\Lambda_{\bR}$, and the subgroup 
$\Gamma$ of $\tLambda$ contains the subgroup generated by the
vertices of $Q$ as a subgroup of finite index. Moreover, 
$\Aut^T(X,L) = \Hom(\Gamma,\bG_m)$. 

Finally, given a finite-dimensional $T$-module $V$, the toric
varieties over $\bP(V)$ correspond to those pairs $(\Gamma,Q)$ where
the vertices of $Q$ are weights of $V^*$, i.e., opposites of weights
of $V$. Equivalently, the moment polytopes of toric varieties over
$\bP(V)$ are exactly the convex hulls of opposites of weights of $V$.

\medskip

\noindent
2) \emph{Reductive varieties}.

Here the acting group is $G \times G$ with Borel subgroup 
$B^- \times B$, maximal torus $T\times T$, weight group 
$\Lambda \times \Lambda$, and set of dominant weights
$(-\Lambda^+) \times \Lambda^+$. The affine reductive varieties of
\cite{AB-I} are the affine spherical $G\times G$-varieties whose
weight group is a direct summand of the anti-diagonal 
$\{(-\lambda,\lambda)~\vert~ \lambda \in \Lambda\} \simeq \Lambda$. 
These varieties are classified by \emph{$W$-admissible cones}, i.e.,
rational polyhedral convex cones $\sigma$ in $\Lambda_{\bR}$ such that 

\begin{enumerate}

\item
the relative interior $\sigma^0$ meets $\Lambda^+_{\bR}$, and 

\item
the distinct $w\sigma^0$ ($w\in W$) are disjoint. 
\end{enumerate}

The weight lattice of the reductive toric variety $X_{\sigma}$ is the
group generated by $\sigma \cap \Lambda$, and its weight cone is 
$\sigma \cap \Lambda^+_{\bR}$. 

Likewise, the polarized reductive varieties of \cite{AB-II} are the
PSV's for $G\times G$ such that their affine cone is a reductive
variety. These varieties are classified by \emph{$W$-admissible polytopes},
i.e., convex lattice polytopes $\delta$ in $\Lambda_{\bR}$ satisfying
the above conditions. The moment polytope of the polarized reductive
variety $(X_{\delta},L_{\delta})$ is the intersection of $\delta$ with
the positive Weyl chamber.

An example of a $W$-admissible polytope is the convex hull of the Weyl 
group orbit of a dominant weight $\lambda$. The vertices of the
corresponding moment polytope are $\lambda$ and its projections on the
faces of the positive chamber. This yields many examples of moment
polytopes with non-integral vertices. Specifically, if $G = \SL(n+1)$
with fundamental weights $\omega_1,\omega_2,\ldots,\omega_n$, and 
$\lambda = \omega_1$, then the vertices are 
$\omega_1, \frac{1}{2}\omega_2,\ldots,\frac{1}{n}\omega_n$.
More complicated examples exist in which the moment polytope of a
spherical variety is not obtained by intersecting a $W$-invariant
integral convex polytope with the positive chamber.

Also, note that the data of the weight lattice and weight cone 
(or moment polytope) do not suffice to distinguish reductive varieties.
This happens, e.g., for $G = \SL(2)$ and the subvarieties of 
$\bP(M_2 \oplus k)$ defined by the equations $ad - bc = z^2$, 
resp.~$ad - bc = 0$. Here $M_2$ denotes the space of $2\times 2$
matrices with coefficients $a,b,c,d$, where $G\times G$ acts by left
and right multiplication; and $k$ denotes the trivial 
$G\times G$-module with coordinate $z$.

\medskip

\noindent
3) \emph{Spherical varieties for $\SL(2)$}.

Let $G = \SL(2)$ with Borel subgroup $B$ of upper triangular matrices,
and maximal torus $T$ of diagonal matrices. We identify $T$ to $\bG_m$
via $t\mapsto \diag(t,t^{-1})$. This identifies $\Lambda$ to $\bZ$,
and $\Lambda^+$ to the set $\bN$ of non-negative integers.
Each simple $G$-module $V(n)$ may be realized as the space
of homogeneous polynomials of degree $n$ in two variables $x$, $y$,
where $G$ acts by linear change of variables.

The classification of the non-trivial projective spherical
$G$-varieties $X$ together with their weight group 
$\Gamma = \Lambda(X)$ is easy; the results are as follows. 

\smallskip

\noindent
(i) $X$ is the projective line $\bP^1 = G/B$. Here $\Gamma$ is trivial.

\smallskip

\noindent
(ii) $X$ is a rational ruled surface
$\bF_e := \bP(\cO_{\bP^1}\oplus \cO_{\bP^1}(e))$ ($e\ge 1$),
where $G$ acts via its natural action on $\bP^1$ and on the
linearized sheaf $\cO_{\bP^1}(e)$. Then $X$ contains a dense
$G$-orbit with isotropy group $U_e$, the semi-direct product of
$U$ with the subgroup of $e$-th roots of unity in $T$. The complement
of this orbit consists of two closed orbits $C_+$, $C_-$, sections of
the fibration $\bF_e \to \bP^1$. The self-intersection of $C_{\pm}$ is 
$\pm e$. Here $\Gamma = e \bZ$.

\smallskip

\noindent
(iii) $X$ is a normal projective surface $S_e$ obtained from $\bF_e$ by
contracting the negative section $C_-$. Then $X$ still contains a
dense orbit with isotropy group $U_e$; its complement is the disjoint
union of a fixed point and the positive section $C_+$. Again, 
$\Gamma = e \bZ$. 

\smallskip

\noindent
(iv) $X = \bP^1\times \bP^1$ where $G$ acts diagonally. The orbits
are the diagonal and its complement, a dense orbit isomorphic to $G/T$.
Here $\Gamma = 2 \bZ$.

\smallskip

\noindent
(v) $X = \bP^2$ where $G$ acts by the projectivization of its linear
action on the quadratic forms in two variables. The orbits are the 
conic of degenerate forms and its complement, a dense orbit isomorphic
to $G/N_G(T)$. Here $\Gamma = 4 \bZ$.

\smallskip

Note that $S_1 \simeq \bP^2$ as abstract varieties, but not as
$G$-varieties (since the orbit structures are different).

We now describe the ample invertible sheaves $L$ on these varieties
$X$ and the corresponding moment polytopes $Q=Q(X,L)$ and weight
groups $\tGamma = \tLambda(X,L)\subseteq \bZ \times \Lambda = \bZ^2$ 
(recall that any invertible sheaf on a normal $G$-variety admits a
unique $G$-linearization, since $G$ is semi-simple and simply
connected).

\smallskip

\noindent
(i) $L = \cO_{\bP^1}(n)$, where $n\ge 1$. Then $Q$ is just the
point $n$. We may realize $X$ as the orbit $G\cdot [x^n]$ in
$\bP(V(n))$; then $L$ is the restriction of $\cO(1)$, and 
$\tGamma = \bZ (1,n)$. 

\smallskip

\noindent
(ii) The degrees of the restriction of $L$ to the sections $C_{\pm}$
are positive integers $n_{\pm}$. One checks that $L$ is uniquely
determined by these integers, which satisfy $n_- < n_+ $ and 
$n^+ - n^-$ is divisible by $e$; we write $L = \cO(n_-,n_+)$. Then the
moment polytope is the interval $[n_-,n_+]$. We may realize $X$ as the
closure of the orbit $G \cdot [x^{n_-} + x^{n_+}]$ in 
$\bP(V(n_-) \oplus V(n_+))$; then $L$ is the restriction of $\cO(1)$.
The group $\tGamma$ is generated by $(1,n_+)$, $(0,e)$.

\smallskip

\noindent
(iii) Likewise, $L$ is uniquely determined by the degree $n$ of its
restriction to $C_+$, a positive integer; we write $L = \cO(n)$. 
Then $Q = [0,n]$. Moreover, $X$ is the closure of the orbit  
$G\cdot [1 + x^n]$ in $\bP(V(0) \oplus V(n))$, and $L$ is the
restriction of $\cO(1)$. The group $\tGamma$ is generated by $(1,n)$,
$(0,e)$.

\smallskip

\noindent
(iv) $L = \cO(m,n)$ with obvious notation, where $m,n\ge 1$. Then
$Q = [\vert m-n\vert, m+n]$. We may realize $X$ as the closure of the
orbit $G\cdot [x^m \otimes y^n]$ in $\bP(V(m)\otimes V(n))$; then $L$ is
the restriction of $\cO(1)$. The group $\tGamma$ is generated by
$(1,m+n)$ and $(0,2)$. 

\smallskip

\noindent
(v) $L = \cO_{\bP^2}(n)$, where $n\ge 1$. Then $Q = [0,2n]$.
We may realize $X$ as the closure of the orbit $G\cdot [x^ny^n]$ in
$\bP(V(2n))$; then $L$ is the restriction of $\cO(1)$. The group 
$\tGamma$ is generated by $(1,2n)$ and $(0,4)$. 

\smallskip

Note that $(\bP^1 \times \bP^1,\cO(m,n))$ may be realized as a
spherical variety over $\bP(V(m+n))$, via the multiplication map
$V(m)\otimes V(n) \to V(m+n)$. The corresponding morphism
$f: \bP^1 \times \bP^1 \to \bP(V(m+n))$
is injective with non-normal image if $m\ne n$, resp. of degree $2$
with smooth image if $m=n$.

In case (iv), both vertices of the moment polytope lie in $\tGamma$,
but the vertex $\vert m-n \vert$ does not arise from any
$G$-subvariety. On the other hand, in case (v) the vertex $0$ is not
in $\tGamma$ for odd $n$.

Also, note that the PSV's $(\bP^1 \times \bP^1,\cO(1,1))$ and
$(S_2,\cO(2))$ have the same weight lattice and moment polytope, but
are non-isomorphic; similarly for $(\bP^2,\cO(2))$ and $(S_4,\cO(4))$.
In both cases, the nonsingular variety can be $G$-equivariantly
degenerated to the singular one.

Finally, given a finite-dimensional $G$-module $V$ with weight set
$F\subset \bN$, the moment polytopes of spherical varieties over
$\bP(V)$ are the intervals with endpoints in $F$ and also the
intervals $[n-2p,n]$, where $n\in F$, $p\in \bN$, and $2p \le n$.

\section{Stable spherical varieties}

\subsection{Affine SSV's}

In this section, we begin by formulating a general definition of
stable spherical varieties. Then we obtain several characterizations
and properties of these varieties, in the affine case. 

We fix a subgroup $\Gamma$ of the weight group $\Lambda$.

\begin{definition}\label{stable}
A \emph{stable spherical variety with weights in $\Gamma$} is a
$G$-variety $X$ satisfying the following conditions :

\smallskip

\noindent
(i) $X$ is semi-normal.

\smallskip

\noindent
(ii) $X$ contains only finitely many $G$-orbits, and these are spherical. 

\smallskip

\noindent
(iii) For any $G$-orbit closure $Y$, the weight group $\Lambda(Y)$ is
contained in $\Gamma$ as a direct summand.

\smallskip

We then say for brevity that $X$ is a SSV.
\end{definition}

\begin{remark}
This definition should be compared with that of a stable toric variety
under a torus $T$ (in the sense of \cite[Definition 1.1.5]{Al}),
i.e., a semi-normal variety where $T$ acts with finitely many orbits
and connected isotropy groups. Our conditions (i) and (ii) are direct
analogs of their toric versions, while (iii) is a more hidden analog
of the connectedness of isotropy groups. This will be developed in
Lemma \ref{stv}.
\end{remark}

By Proposition \ref{psv}, any spherical variety is a SSV. Conversely,
any SSV is obtained by gluing spherical $G$-varieties along
$G$-subvarieties:

\begin{proposition}\label{ssv}
Let $X$ be a SSV and denote by $\cQ(X)$ the finite set of $G$-orbit
closures in $X$, partially ordered by inclusion. Then every such orbit
closure $Y$ is a spherical variety. Moreover, the natural map
$$
f: \varinjlim_{Y \in \cQ(X)} Y \to X
$$
is an isomorphism.
\end{proposition}

\begin{proof}
Let $\nu: Y' \to Y$ be the normalization, $\cO$ a $G$-orbit in
$Y$, and $Z:=\bar{\cO}$ its closure. Put $Z' := \nu^{-1}(Z)$, then
$Z'$ is irreducible by Lemma \ref{cover}. Further,
$\Lambda(Z)\subseteq \Lambda(Z')$ as a subgroup of finite index,
by Lemma \ref{weightset} (ii). On the other hand, $\Lambda(Z)$ is a
direct summand of $\Gamma$, and Lemma \ref{weightset} (i) yields
$$
\Lambda(Z') \subseteq \Lambda(Y')=\Lambda(Y) \subseteq \Gamma.
$$
It follows that $\Lambda(Z') = \Lambda(Z)$. By Lemma
\ref{cover} again, this implies that $\nu$ restricts to an
isomorphism $\nu^{-1}(\cO) \to \cO$. Since this holds for all orbits,
$\nu$ is bijective. Moreover, since $Z'$ is normal by Proposition
\ref{psv}, the restriction $\nu\vert_{Z'}$ is the normalization. This
yields a map 
$$
f': \varinjlim_{Y' \in \cQ(X)} Y' \to X,
$$
which is clearly finite and bijective. (Since the morphisms 
$Z' \to Y'$ are closed immersions, the limits in the categories of
schemes and of affine schemes coincide). As $X$ is semi-normal, $f'$
is an isomorphism; this implies both statements.
\end{proof}

\begin{remark}
One shows by similar arguments that the $G$-varieties obtained as
direct limits of a finite poset of spherical varieties (with arrows
being closed immersions) are characterized by the conditions (i) and
(ii) of Definition \ref{ssv}, together with

\smallskip

\noindent
(iii)' For any $G$-orbit closures $Z\subseteq Y$, the weight group
$\Lambda(Z)$ is a direct summand of $\Lambda(Y)$.

\smallskip

However, there exist $G$-varieties which satisfy (i), (ii), (iii)',
but such that (iii) fails for any subgroup $\Gamma$ of $\Lambda$ (for
example, take $G=\bG_m$ and let $X$ be the union of two affine lines
where $G$ acts linearly with weights $2$, $3$). 
\end{remark}

Until the end of this section, we will only consider \emph{affine}
SSV's. We begin by relating them to stable toric varieties. For this,
recall that each affine $G$-variety $X=\Spec \, R$ defines an affine
$T$-variety
$$
X//U := \Spec \, R^U,
$$ 
where $R^U$ denotes the ($T$-stable) subalgebra of $U$-invariants in
$R$. Clearly, $X$ and $X//U$ have the same weight set. In fact, many
properties of $X$ can be read off $X//U$. This is developed in
\cite[\S 6]{Po} and \cite[Chapter 3]{Gr}; we will freely use the
results there.

\begin{lemma}\label{stv}
The following conditions are equivalent for an affine $G$-variety $X$:

\smallskip

\noindent
(i) $X$ is a SSV with weights in $\Gamma$.

\smallskip

\noindent
(ii) The action of $T = \Hom(\Lambda,\bG_m)$ on $X//U$ factors through
an action of the quotient torus $T_{\Gamma} := \Hom(\Gamma,\bG_m)$ for
which $X//U$ is a stable toric variety. 
\end{lemma}

\begin{proof}
(i) $\Rightarrow$ (ii) By \cite[Lemma 2.1]{AB-I}, $X//U$ is a
semi-normal $T$-variety. Its irreducible components are the $Y//U$,
where $Y$ runs over the irreducible components of $X$. Since $Y$ is
normal by Proposition \ref{ssv}, each $Y//U$ is an affine toric
variety. And since the weight group $\Lambda(Y)$ is a direct summand
of $\Gamma$, the $T$-action on $Y//U$ factors through a
$T_{\Gamma}$-action having a dense orbit with connected isotropy
group: the subtorus of $T_{\Gamma}$ with character group
$\Gamma/\Lambda(Y)$. It follows that any $T_{\Gamma}$-isotropy group
in $Y//U$ is connected. This means that $X//U$ is a
$T_{\Gamma}$-stable toric variety. 

(ii) $\Rightarrow$ (i) By \cite[Lemma 2.1]{AB-I} again, $X$ is
semi-normal. Moreover, by \cite[Section 2.3]{Al}, $Y//U$ is an affine
$T_{\Gamma}$-toric variety for any irreducible component $Y$ of
$X$. Thus, $Y$ is a spherical $G$-variety. In particular, every
$G$-orbit $\cO$ in $X$ is spherical. Let $Z:=\bar{\cO}$, then  
$\Lambda(Z) = \Lambda(Z//U)$, and $Z//U$ is a $T$-orbit closure in
$X//U$. Since $T$ acts on $Z//U$ via $T_{\Gamma}$ and the
$T_{\Gamma}$-isotropy groups are connected, it follows that
$\Lambda(Z)$ is a direct summand of~$\Gamma$.
\end{proof}

Clearly, any affine SSV is \emph{multiplicity-bounded}, i.e., 
the multiplicities of simple $G$-modules in its affine ring are
uniformly bounded. Those affine SSV's that are multiplicity-free 
admit a useful characterization, which follows from 
Lemma \ref{stv} and \cite[Lemma 2.3]{AB-I}.

\begin{lemma}\label{multiplicityfree}
For an affine multiplicity-free variety $X$, the following
conditions are equivalent:

\smallskip

\noindent
(i) $X$ is a SSV with weights in $\Gamma$.

\smallskip

\noindent
(ii) The weight set $\Lambda^+(X)$ is saturated in $\Gamma$, i.e., 
$\Lambda^+(X) = C(X) \cap \Gamma$.
\end{lemma}

From this lemma, we now deduce important invariance properties of
$G$-subvarieties of SSV's:

\begin{lemma}\label{reduced}
Let $X$ be an affine multiplicity-free SSV. Then:

\smallskip

\noindent
(i) Any $G$-subvariety $Y$ is a multiplicity-free SSV, invariant under
$\Aut^G(X)$.

\smallskip

\noindent
(ii) For any two $G$-subvarieties $Y$, $Z$, the scheme-theoretic
intersection $Y\cap Z$ is non-empty, connected and reduced.
\end{lemma}

\begin{proof}
By Lemma \ref{stv}, we may assume that $G = T = T_{\Gamma}$. Then 
$\Gamma = \Lambda$, and $X$ is a stable toric variety.

(i) By Lemma \ref{multiplicityfree}, to show that $Y$ is a SSV, it
suffices to check that $C(Y) \cap \Lambda \subseteq \Lambda^+(Y)$.
Given $\lambda \in C(Y)\cap \Lambda$, we may find a positive integer
$n$ such that $n\lambda \in \Lambda^+(Y)$. Thus, 
$\lambda \in \Lambda^+(X)$ since the latter is saturated in
$\Gamma$. Let $R$ (resp. $S$) be the affine ring of $X$ (resp. $Y$),
and $I_Y$ the ideal of $Y$ in $R$. By the exact sequence of
$G$-modules 
$$
0 \to I_Y  \to R \to S \to 0
$$
and multiplicity-freeness, $\Lambda ^+(X) = \Lambda^+(R)$ is the
disjoint union of $\Lambda^+(Y)=\Lambda^+(S)$ and $\Lambda^+(I_Y)$. 
Since $\lambda \notin \Lambda^+(I_Y)$ (as 
$n\lambda \notin \Lambda^+(I_Y)$), it follows that 
$\lambda \in \Lambda^+(Y)$ as desired.

To show that $Y$ is invariant under $G$-automorphisms of $X$, we note
that the irreducible components of $X$ are affine toric varieties with
pairwise distinct weight cones (by multiplicity-freeness again): these
components are pairwise non-isomorphic. Thus, each component is
invariant under $\Aut^G(X)$. Together with Lemma \ref{automorphisms},
this completes the proof.

(ii) Every $G$-invariant regular function on $X$ is constant, since
$X$ is multiplicity-free. Thus, $X$ contains a unique closed
$G$-orbit. It follows that $Y \cap Z$ is non-empty and connected.

To show the reducedness, consider a $T$-eigenvector
$f\in R_{\lambda}$ and a positive integer $n$ such that 
$f^n\in I_{Y\cap Z} = I_Y + I_Z$. Since $f^n\in R_{n\lambda}$ and this
weight space is a line, then $f^n \in I_Y$ or $f^n \in I_Z$. Thus,
$f\in I_Y$ or $f\in I_Z$, so that $f\in I_{Y\cap Z}$.
\end{proof}

Finally, we obtain a simple sufficient condition for a stable
spherical variety to be Cohen--Macaulay.

\begin{proposition}\label{CM}
Let $X$ be an affine multiplicity-free SSV. If the weight cone $C(X)$
is convex, then $X$ is Cohen--Macaulay.
\end{proposition}

\begin{proof}
We adapt the argument of \cite[Lemma 5.15]{AB-II}.
By \cite{Po} (see also \cite{Gr}, \cite[Section 2]{AB-III}), $X$
admits a flat degeneration to an affine $G$-variety
$X_0$ which is ``horospherical'', i.e., $X_0 = G \cdot X_0^U$; the
base of this degeneration may be taken to be the affine line.
Moreover, the affine rings of $X$ and $X_0$ are 
isomorphic as $G$-modules, so that $X_0$ is a  multiplicity-free SSV
as well. Since the property of being Cohen--Macaulay is open for
fibers of this degeneration, we may assume that $X$ itself is
horospherical.

Now let $Y:=X^{U^-}$, so that $X = G \cdot Y$. 
By \cite[Lemma 2.4]{AB-III}, the composed map $Y \to X \to X//U$ is a
$T$-isomorphism. Together with Lemma \ref{stv}, it follows that $Y$ is
an affine stable toric variety for $T_{\Gamma}$, and 
$$
\Lambda^+(Y) = \Lambda^+(X) = \Gamma \cap C(X).
$$ 
Thus, $Y$ is Cohen-Macaulay by \cite[Theorem 2.3.19]{Al}. Moreover,
the dualizing module $H^0(Y,\omega_Y)$, regarded as a $T$-module, has
weight set $\Gamma \cap C(X)^0$, where $C(X)^0$ denotes the relative
interior of the cone $C(X)$ (see \cite[p.~405--406]{AB-II}). 

Let $P$ be the set of all $g\in G$ such that $g \cdot Y = Y$; this is
a parabolic subgroup of $G$ containing $B^-$. Moreover, any 
$\lambda \in \Lambda^+(X)$ extends to a character of $P$; thus, the
$B^-$-action on $Y$ (via the quotient $B^- \to B^-/U^- = T$) extends
to an action of $P$. In fact, $P$ is the largest parabolic subgroup
containing $B^-$ such that every $\lambda \in \Lambda^+(X)$ extends to a
character of $P$. Thus, any such $\lambda$ yields a $G$-linearized
invertible sheaf $\cL(\lambda)$ on $G/P$, which is globally
generated. Moreover, $\cL(\lambda)$ is ample for any 
$\lambda \in \Gamma \cap C(X)^0$. 

Let $X'$ denote the fiber product $G\times^P Y$. In 
other words, $X'$ is a $G$-variety equipped with a $G$-morphism
$$
p:X' \to G/P
$$
such that the fiber at $P$ is the $P$-variety $Y$. Then $p$ is a
locally trivial fibration; thus, $X'$ is Cohen--Macaulay. We also
have a proper, surjective $G$-morphism
$$
\pi: X' \to X, \quad (g,y) P \mapsto g\cdot y.
$$
By a lemma of Kempf, to show that $X$ is Cohen--Macaulay, it suffices
to check that
$\pi_*(\cO_{X'}) = \cO_X$, $R^i\pi_*(\cO_{X'}) =0$ for all $i\ge 1$,
and $R^i\pi_*(\omega_{X'}) =0$ for all $i\ge 1$, where $\omega_{X'}$
denotes the dualizing sheaf of $X'$. Since $X$ is affine, 
this is equivalent to 
$H^0(X',\cO_{X'}) = H^0(X,\cO_X)$, 
$H^i(X',\cO_{X'}) = 0$ for all $i\ge 1$, and $H^i(X',\omega_{X'})=0$
for all $i\ge 1$.

Since the morphism $p$ is affine and 
$p_*\cO_{X'} \simeq \bigoplus_{\lambda\in \Gamma \cap C(X)} \cL(\lambda)$,
we obtain 
$$
H^i(X',\cO_{X'}) \simeq 
\bigoplus_{\lambda\in \Gamma \cap C(X)} H^i(G/P,\cL(\lambda))
\text{ for all }i\ge 0.
$$
Thus, the desired assertions on  $H^i(X',\cO_{X'})$ follow from the
Borel-Weil theorem. To compute $H^i(X',\omega_{X'})$, we use the
isomorphism of $G$-linearized sheaves 
$\omega_{X'} \simeq (p^*\omega_{G/P}) \otimes \omega_p$, 
where $\omega_p$ denotes the relative dualizing sheaf. Thus, 
$p_*\omega_{X'} \simeq \omega_{G/P} \otimes p_* \omega_p$.
Moreover, $p_*\omega_p$ is the $G$-linearized sheaf on $G/P$
associated with the $P$-module 
$H^0(Y,\omega_Y) \simeq 
\bigoplus_{\lambda\in \Gamma \cap C(X)^0} \cL(\lambda)$. Thus, 
$$
H^i(X',\omega_{X'}) \simeq 
\bigoplus_{\lambda\in \Gamma \cap C(X)^0} 
H^i(G/P,\omega_{G/P}\otimes\cL(\lambda))
$$
and the latter vanishes for all $i\ge 1$, by the Kodaira vanishing
theorem.
\end{proof}

\begin{remark}
By this argument, an affine multiplicity-free SSV is
Cohen--Macaulay whenever its weight cone $C$ satisfies the following
conditions: 

\begin{enumerate}

\item 
$C$ is homeomorphic to a convex cone.

\item 
$C$ is a finite union of convex cones $C_i$ such that the smallest
face of the positive chamber containing $C_i$ is independent of $i$.

\end{enumerate}

\end{remark}

\subsection{Polarized SSV's} 

In this section, we adapt the definition of stable spherical varieties
to the setting of polarized varieties, and we generalize results of
Sections 1.3 and 2.1 to these polarized SSV's.

We fix a subgroup $\Gamma$ of $\tLambda$.

\begin{definition}
A \emph{polarized stable spherical variety with weights in $\Gamma$}
is a polarized $G$-variety $(X,L)$ satisfying the following
conditions:

\smallskip

\noindent
(i) $X$ is semi-normal.

\smallskip

\noindent
(ii) $X$ contains only finitely many $G$-orbits, and these are
spherical.

\smallskip

\noindent
(iii) $\tLambda(Y,L)$ is contained in $\Gamma$ as a direct summand,
for any $G$-orbit closure~$Y$.

\smallskip

We then say for brevity that $(X,L)$ is a PSSV.
\end{definition}

Now \cite[Lemma 2.1]{AB-II} and Lemma \ref{ssv} imply readily the
following:

\begin{lemma}\label{affinecone}
A polarized $G$-variety $(X,L)$ is a PSSV for $G$ if and only if its 
affine cone is a SSV for $\tG$. Then 
$$
(X,L) = \varinjlim_{Y\in \cQ(X)} (Y,L),
$$ 
where $\cQ(X)$ denotes the poset of $G$-orbit closures in $X$.
\end{lemma}

We may now generalize the results of Propositions \ref{psv},
\ref{reduced} and \ref{CM} to the polarized setting. From now on we
only consider \emph{multiplicity-free} SSV's. Indeed, we are mainly 
interested in spherical varieties and their stable limits, which are
all multiplicity-free. The general SSVs can be studied along the
lines of \cite[Section 2]{Al}.

\begin{proposition}\label{pssv}
Let $(X,L)$ be a PSSV with weights in $\Gamma$ and moment set $Q$,
and let $Y$ be a $G$-subvariety of $X$. Then: 

\smallskip

\noindent
(i) $(Y,L)$ is a PSSV, invariant under $\Aut^G(X,L)$.

\smallskip

\noindent
(ii) The scheme-theoretic intersection $Y\cap Z$ is reduced, for any
$G$-subvariety $Z$ of $X$.  

\smallskip

\noindent
(iii) The restriction map $H^0(X,L^n) \to H^0(Y,L^n)$ is surjective
for all $n\ge 0$.

\smallskip

\noindent
(iv) As $G$-modules, 
$$
H^0(X,L^n) \simeq \bigoplus_{\lambda\in\Lambda^+,\;
(n,\lambda)\in \Gamma \cap nQ} V(\lambda).
$$ 
Equivalently, as $\tG$-modules,
$$
R(X,L) \simeq \bigoplus_{\tlambda \in \Gamma \cap \Cone(Q)}
V(\tlambda).
$$

\smallskip

\noindent
(v) $L$ is globally generated.

\smallskip

\noindent
(vi) $H^i(X,L^n) = 0$ for all $i,n\ge 1$.

\smallskip

\noindent
(vii) If $Q$ is convex, then the affine cone over $(X,L)$ is
Cohen--Macaulay. In particular, $X$ is Cohen--Macaulay.
\end{proposition}

\begin{proof}
(i) and (ii) follow from Lemmas \ref{reduced} and \ref{affinecone}.

(iii) Consider the finite injective map $f:\tY \to \tX$. By Lemma
\ref{reduced} again, $f(\tY)$ is semi-normal. Thus, $f$ is a closed
immersion. In other words, the restriction map $R(X,L) \to R(Y,L)$ is
surjective. 

(iv) is a consequence of Lemma \ref{multiplicityfree}.  

(v) follows from (iii) as in the proof of Lemma \ref{psv}.

(vi) Let $Y$ be an irreducible component of $X$, and $Z$ the union of
all the other irreducible components. Consider the Mayer--Vietoris
exact sequence
$$
0 \to \cO_X \to \cO_Y \oplus \cO_Z \to \cO_{Y\cap Z} \to 0.
$$
Note that $(Y,L)$ is a PSV, and the connected components of $Z$,
$Y\cap Z$ (equipped with the restriction of $L$) are multiplicity-free
PSSV's, strictly smaller than $X$. By Lemma \ref{psv} and induction,
this yields an exact sequence
$$\displaylines{
0 \to H^0(X,L^n) \to H^0(Y,L^n) \oplus H^0(Z,L^n) \to H^0(Y\cap Z,L^n)
\hfill\cr\hfill
\to H^1(X,L^n) \to 0
\cr}
$$
and also $H^i(X,L^n) = 0$ for all $i\ge 2$, $n\ge 1$. To complete the 
proof, it suffices to check that the restriction map
$$
H^0(Y,L^n) \to H^0(Y\cap Z,L^n)
$$
is surjective for $n\ge 1$. Since $Y\cap Z$ may be disconnected, the
statement (iii) does not apply readily. We circumvent this by
introducing
$$
R(Y\cap Z, L) := k \oplus \bigoplus_{n=1}^{\infty} H^0(Y\cap Z,L^n)
$$
(a reduced, finitely generated algebra) and 
$$
\widetilde{Y\cap Z}:= \Spec \, R(Y\cap Z,L).
$$ 
Then we obtain a finite injective map 
$f: \widetilde{Y\cap Z} \to \tY$, with image a (connected)
$\tG$-subvariety. By Lemma \ref{reduced} again, $f$ is a closed
immersion; this is equivalent to the desired surjectivity.

(vii) follows from Proposition \ref{CM}.
\end{proof}

\subsection{Structure}

In this section, we obtain a partial generalization of the
classifications of stable toric varieties (in \cite{Al}) and of stable 
reductive varieties (in \cite{AB-I,AB-II}), to the setting of stable
spherical varieties. We only treat the polarized case, the affine case
being entirely similar, and we still assume multiplicity-freeness
throughout. The classification of SSV's over a fixed projective space
will be given in Lemma \ref{maps}.

\begin{lemma}\label{subdivision}
Let $(X,L)$ be a PSSV with poset $\cQ$ of $G$-orbit closures, and
moment set $Q$. Then: 

\smallskip

\noindent
(i) $Q = \bigcup_{Y \in \cQ} Q(Y,L)$.

\smallskip

\noindent
(ii) $Y\cap Z \in \cQ$ for any $Y$, $Z$ in $\cQ$. Moreover, 
$Q(Y\cap Z,L)$ is the unique common face to $Q(Y,L)$ and $Q(Z,L)$.

\smallskip

\noindent
(iii) $Z\subseteq Y$ if and only if $Q(Z,L)$ is a face of $Q(Y,L)$. In 
particular, $Q(Y,L)$ determines $Y$ uniquely.
\end{lemma}

\begin{proof}
Recall that $Q$ is the union of the moment polytopes of the 
irreducible components of $X$. Moreover, $Q(Z,L)$ is a face of
$Q(Y,L)$ whenever $Z\subseteq Y$. This shows (i) and one implication
of (iii). For the converse implication, if $Q(Z,L)\subseteq Q(Y,L)$,
then the same inclusion holds for the weight sets of the affine cones
$\tZ$, $\tY$. By multiplicity-freeness, the opposite inclusion holds
for the ideals of $\tZ$, $\tY$ in $R(X,L)$. Thus, $\tZ \subseteq \tY$
and $Z \subseteq Y$.

(ii) By (the proof of) Proposition \ref{pssv} (vi), the restriction
maps 
$$
R(Y,L) \to R(Y\cap Z,L), \quad R(Z,L) \to R(Y\cap Z,L)
$$ 
are both surjective. Thus, the Mayer--Vietoris sequence
$$
0 \to R(Y\cup Z, L) \to R(Y,L) \oplus R(Z,L) \to R(Y\cap Z,L) \to 0
$$ 
is exact. Moreover, $R(Y\cup Z, L)$, $R(Y,L)$, $R(Z,L)$, $R(Y\cap Z,L)$
are all multiplicity-free $\tG$-modules, and the corresponding weight
sets are saturated in $\Gamma$. It follows that 
$\tC(Y\cap Z,L) = \tC(Y,L) \cap \tC(Z,L)$, so that 
$$
Q(Y\cap Z,L) = Q(Y,L) \cap Q(Z,L).
$$ 
In particular, $Q(Y\cap Z, L)$ is convex. On the other hand, 
$Q(Y\cap Z, L)$ is a union of faces of $Q(Y,L)$ corresponding to the
irreducible components of $Y\cap Z$, by Proposition \ref{psv}. Thus,
$Q(Y\cap Z,L)$ is a unique face of $Q(Y,L)$, and $Y\cap Z$ is
irreducible.
\end{proof}

With the preceding assumptions and notation, the moment polytopes
$Q(Y,L)$, $Y\in \cQ$, form a partial subdivision of $Q$ by rational
convex polytopes. The characterization of those subdivisions that arise
from PSSV's is an open problem; some examples are considered in 2.4. 

We now describe those PSSV's $(X,L)$ that have prescribed building
blocks $(Y,L)$, where $Y\in \cQ$. Recall from Proposition \ref{ssv}
and Lemma \ref{affinecone} that
$$
(X,L) \simeq \varinjlim_{Y \in \cQ} (Y,L_Y),
$$
where $L_Y$ denotes the restriction of $L$ to $Y$, and 
the directed system is defined by the inclusions 
$$
i_{ZY}:(Z,L_Z) \to (Y,L_Y) \quad (Z\subset Y).
$$
We introduce twists of this directed system, as follows. Let
$t=(t_{ZY})_{Z\subset Y}$, where the $t_{ZY} \in \Aut^G(Z,L_Z)$ satisfy
$$
t_{Y_3Y_1} = t_{Y_2Y_1}\vert_{Y_3} \circ t_{Y_3Y_2}
$$
whenever $Y_3 \subset Y_2 \subset Y_1$ 
(recall that the group $\Aut^G(Y_2,L_{Y_2})$ is abelian and 
leaves $(Y_3,L_{Y_3})$ invariant, see Lemma \ref{automorphisms}). 
This defines a twisted direct system
$$
i_{ZY} \circ t_{ZY} : (Z,L_Z) \to (Y,L_Y).
$$
We denote its direct limit by $(X_t,L_t)$.

Note that the twists $t$ are precisely the $1$-cocycles of the complex
of abelian groups $C^*(\cQ,\Aut)$:
$$ \prod_{Y_1} \Aut^G(Y_1,L_{Y_1}) \to 
\prod_{Y_2 \subset Y_1} \Aut^G(Y_2,L_{Y_2}) \to 
\prod_{Y_3 \subset Y_2 \subset Y_1} \Aut^G(Y_3,L_{Y_3}) 
\to \cdots
$$
with the usual differentials. This easy yields

\begin{lemma}\label{cohomology}
For any $1$-cocycle $t\in Z^1(\cQ,\Aut)$, the limit $(X_t,L_t)$ is a
PSSV having the same building blocks $(Y,L_Y)$,
$Y\in \cQ$, as $(X,L)$. Moreover, $\Aut^G(X_t,L_t) = H^0(\cQ,\Aut)$,
and the isomorphism classes of the polarized $G$-varieties $(X_t,L_t)$
are parametrized by the first cohomology group $H^1(\cQ,\Aut)$
(a diagonalizable linear algebraic group).
\end{lemma}

\begin{remark}
The complex $C^*(\cQ,\Aut)$ admits an interpretation in terms of
sheaf cohomology on the poset $\cQ$. Indeed, consider $\cQ$ as a
topological space, where the open subsets are the decreasing subsets,
i.e., the unions of subsets 
$$
\cQ_{\le Y} := \{ Z \in \cQ ~\vert~ Z \subseteq Y\}.
$$
Then a sheaf $\cF$ of abelian groups on $\cQ$ is given by a collection
of abelian groups $\cF(Y)$, where $Y\in \cQ$ (the stalks of $\cF$),
together with restriction maps $\cF(Y) \to \cF(Z)$, where $Z\subset Y$,
which are compatible in an obvious sense. This yields a complex
$C^*(\cQ,\cF)$ which sheafifies to a resolution of $\cF$ by flabby
sheaves, analogous to Godement's canonical resolution. Thus, the 
$H^i(\cQ,\Aut)$ are the cohomology groups of the sheaf 
$\Aut: Y \mapsto \Aut^G(Y,L_Y)$.

Observe that $H^i(\cQ_{\le Y}, \cF) =0$ for any sheaf $\cF$ and any 
$i\ge 1$ (since the functor 
$\cF \mapsto H^0(\cQ_{\le Y}, \cF) = \cF(Y)$ 
is exact). Moreover, by Lemma \ref{subdivision} (ii), the intersection of
any family of basic open subsets $\cQ_{\le Y}$ is a basic open
subset. Thus, the cohomology groups of $\cF$ are those of the Cech
complex associated with the covering $\cQ = \bigcup_i \cQ_{\le Y_i}$,
where the $Y_i$ are the irreducible components of $X$:
$$
\bigoplus_i \cF(Y_i) \to \bigoplus_{i<j} \cF(Y_{ij}) \to
\bigoplus_{i<j<k} \cF(Y_{ijk}) \to \cdots
$$
where the $Y_{ij} = Y_i \cap Y_j$ are the double intersections,
$Y_{ijk} = Y_i \cap Y_j \cap Y_k$ are the triples intersections, etc.
\end{remark}

Also, note that the group $H^1(\cQ,\Aut)$ may be infinite. In fact, this
already happens for stable toric surfaces, see 
\cite[Example 2.2.10]{Al} and also Example 2.4(1). 

Together with Lemma \ref{class},  it  follows that there are
infinitely many orbits of projective stable toric surfaces over some
$\bP(V)$, for the action of $\Aut^G \bP(V)$. In other words, Theorem
\ref{finiteness} does not extend to all stable spherical varieties;
however, it holds for SSV's which contain a unique closed $G$-orbit
(see Example 2.4(5)). 

For arbitrary SSV's, we obtain a weaker finiteness statement: 

\begin{lemma}\label{finitefibers}
Let $X$ be a $G$-subvariety of $\bP(V)$, where $V$ is a
finite-dimensional $G$-module. Then there are only finitely many
SSV's $(X',f)$ over $\bP(V)$ such that $f(X') = X$. 
\end{lemma}

\begin{proof}
Let $Y'$ be a $G$-orbit closure in $X'$, with image $Y$ in $X$. 
Then there are only finitely many possibilities for
$(Y',f\vert_{Y'})$, by Lemma \ref{cover}. Moreover, 
the assignment $Y'\mapsto Y$ is a bijection from the poset of 
$G$-orbit closures in $X'$ to the poset $\cQ$ of $G$-orbit closures
in $X$, by Lemmas \ref{cover} and \ref{subdivision}. So we may fix the
building blocks $(Y',L'_{Y'})$ of the pairs $(X',L')$ over
$(X,L)$, where $L := \cO(1)$ and $L' := f^*L$.

Fix such a pair $(X',L')$, so that an arbitrary pair may be written as
$(X'_t,L'_t)$, where $t\in Z^1(\cQ,\Aut)$. Let $t = (t_{Z'Y'})$, then 
$$
t_{Z'Y'}\in \Aut^G(Z',L') \subseteq \Hom(\tLambda(Z',L'),\bG_m).
$$ 
Moreover, the restriction of $t_{Z'Y'}$ to the subgroup
$\tLambda(Z,L)$ of $\tLambda(Z',L')$ is constant, since 
$(X',L')$ and $(X'_t,L'_t)$ map to the same pair $(X,L)$. Thus,
$t_{Z'Y'}$ belongs to a finite subgroup of $\Aut^G(Z',L')$. So there
are only finitely many possibilities for $t$.
\end{proof}

\subsection{Examples}

We illustrate the notions and results of Section 2 in the cases
considered in 1.4, and in two additional cases.
 
\medskip

\noindent
1) \emph{Stable toric varieties}. 

Let $G=T$, then the stable spherical varieties with weights in
$\Lambda$ are exactly the stable toric varieties of \cite{Al}
(for the quotient $T_{\Gamma}$ of $T$), as follows from Lemma
\ref{stv}. The moment polytopes of orbit closures of a
multiplicity-free polarized stable toric variety form a subdivision of
the moment set $Q$ into convex lattice polytopes. Any such subdivision
arises from a multiplicity-free polarized stable toric variety, by
\cite[Section 2.4]{Al}. 

In the case where $\cQ$ is the subdivision of a triangle given in 
\cite[Example 2.2.10]{Al}, one obtains $H^1(\cQ,\Aut) \simeq \bG_m$.
This yields infinitely many stable toric surfaces which share the same
moment set, but are pairwise non-isomorphic.

\medskip

\noindent
2) \emph{Stable reductive varieties}.

With the notation of Example 1.4(2), the stable reductive varieties 
are exactly the SSV's for $G\times G$ with weights in the
anti-diagonal. They are classified by complexes of $W$-admissible
cones in the affine case, resp. polytopes in the polarized case; see
\cite[Section 5]{AB-I} and \cite[Section 2]{AB-II}. Again, there are
non-isomorphic multiplicity-free stable spherical varieties having the
same weight set and the same subdivision of their moment polytope.

\medskip

\noindent
3) \emph{PSSV's for $\SL(2)$}.

We consider projective multiplicity-free SSV's for the group 
$G = \SL(2)$. By Example 1.4(3) and Lemma \ref{subdivision}, these
varieties are exactly the chains 
$$
X = X_1 \cup X_2 \cup \cdots \cup X_r
$$ 
of spherical $G$-varieties, such that: 

\smallskip

\noindent
$X_2,\ldots,X_r$ are all isomorphic to the same rational ruled
surface $\bF_e$,

\smallskip

\noindent
$X_1$ is isomorphic to $\bF_e$, or to $S_e$, or to 
$\bP^1 \times \bP^1$ (if $e=2$), or to $\bP^2$ (if $e=4$),

\smallskip

\noindent
for $i=2,\ldots, r-1$, the positive section of $X_i$ is identified to
the negative section of $X_{i+1}$, 

\smallskip

\noindent
likewise, the positive section of $X_1$ (in the case of $\bF_e$ or
$S_e$) or its closed orbit (in the other cases) is identified to the
negative section of $X_2$.

\smallskip

(Such identifications are unique, since the $G$-automorphism group of
the projective line is trivial). 

Moreover, any ample, $G$-linearized line bundle $L$ on $X$ is
uniquely determined by the degrees of its restrictions to the various
sections. These degrees form a strictly increasing sequence 
$(n_0,n_1,\ldots,n_r)$ in the case of $\bF_e$ or $S_e$, resp. 
$(n_1,\ldots,n_r)$ in the other cases. The moment set $Q=Q(X,L)$ is
the concatenation of the moment intervals of the $X_i$. So 
$Q = [n_0,n_r]$ in the case of $\bF_e$ or $S_e$; $Q=[n_1-2p,n_r]$ for
some integer $p$ such that $0\le p \le \frac{n_1}{2}$ in the case of
$\bP^1 \times \bP^1$; and $Q = [0,n_r]$ in the case of $\bP^2$.

\medskip

\noindent
4) \emph{Some PSSV's for $\SL(2) \times \SL(2)$.}

Let $G = \SL(2)\times SL(2)$. We construct a multiplicity-free PSSV
$(X,f)$ over some projective space $\bP(V)$, such that $X$ is the
union of two spherical varieties $X_1$, $X_2$ glued along a spherical
subvariety $X_{12}$; the corresponding moment polytopes $Q_1$, $Q_2$,
$Q_{12}$ form a subdivision of a convex polytope $Q$, but there exists
no spherical variety over $\bP(V)$ with moment polytope $Q$. 

For this, we adapt an unpublished example of D.~Luna who 
showed that the moduli space of affine spherical varieties with a
prescribed weight monoid (defined in \cite[Section 1.4]{AB-III}) may
have several irreducible components.

The weight group of $G$ is $\bZ^2$ and the subset of dominant weights
is $\bN^2$. Each simple $G$-module $V(n,n')$ is the space
of polynomials in the variables $x,y,x',y'$ which are homogeneous of
degree $n$ in $x,y$, and of degree $n'$ in $x',y'$; then $x^n x^{'n'}$
is a highest weight vector. Consider the $G$-module
$$
V := V(2,0) \oplus V(4,2).
$$
Then $\bP(V)$ has two closed $G$-orbits, $G\cdot [x^2]$ and 
$G\cdot [x^4 x^{'2}]$. 

Let $x_{12}:= [x^2 + x^4 x^{'2}]$ and 
$X_{12}:=\overline{G\cdot x_{12}}$. Then $X_{12}$ is a spherical
subvariety of $\bP(V)$. Let $L_{12}$ be the restriction of
$\cO(1)$ to $X_{12}$, then the weight group 
$\tLambda(X_{12},L_{12})\subset \bZ^3$ is generated by $(1,2,0)$ and
$(1,4,2)$. Moreover, the moment polytope $Q_{12}$ is the segment
$[(2,0),(4,2)]$. Any $G$-subvariety of $\bP(V)$ which contains both
closed orbits must also contain $X_{12}$. 

Next let $x_1 := [xy + x^4 x^{'2}]\in \bP(V)$ and let $X_1$ be the
normalization of $\overline{G\cdot x_1}$ (one can show that the latter
orbit closure is non-normal). Then the $G$-isotropy 
group of $x_1$ is the semi-direct product of the additive group 
$(x\mapsto x, y \mapsto y, x'\mapsto x', y'\mapsto y' + u x')$
with the diagonalizable group 
$(x \mapsto tx, y\mapsto t^{-1}y, x' \mapsto \varepsilon t^{-2}x',
y' \mapsto \varepsilon t^2 y')$, 
where $\varepsilon^2 = 1$. It follows that $X_1$ is a spherical variety. 

Let $L_1$ be the pull-back of $\cO(1)$ to $X_1$. Then, by
restricting to the open orbit $G\cdot x_1$, 
one checks that the weight group $\tLambda(X_1,L_1)$ is 
generated by the three vectors $(1,0,0)$, $(1,2,0)$, $(1,4,2)$;
moreover, the weight cone $\tC(X_1,L_1)$ is contained in the cone
spanned by these vectors. It follows that the moment polytope 
$Q_1 := Q(X_1,L_1) \subset \bZ^2$
is the triangle with vertices $(0,0)$, $(2,0)$, $(4,2)$. Indeed,
we just saw that $Q_1$ is contained in this triangle. On the other
hand, $\overline{G \cdot x_1}$ contains the $G$-subvarieties 
$X_{12}$ and $Y_1 :=\overline{G\cdot [xy]}$. Thus, 
$Q_1$ contains as faces the moment polytopes of $X_{12}$ and $Y_1$,
i.e., the segments $[(2,0),(4,2)]$ and $[(0,0),(2,0)]$.

Next let $x_2 := [x^2 + x^4 x'y']\in \bP(V)$ and let $H_2\subseteq G$
be the subgroup 
$(x\mapsto \varepsilon x, y\mapsto \varepsilon y + u x, 
x'\mapsto tx', y' \mapsto t^{-1}y')$,
where $\varepsilon^2 = 1$. Then $H_2$ is a subgroup of index $2$ of
the $G$-isotropy group of $x_2$. Denote by $X_2$ the normalization of 
$\overline{G\cdot x_2}$ in the function field $k(G/H_2)$, then $X_2$
is a spherical variety.

Let $L_2$ be the pull-back of $\cO(1)$ to $X_2$. One checks
that $\tLambda(X_2,L_2) = \tLambda(X_1,L_1)$. Moreover, the moment
polytope
$Q_2 := Q(X_2,L_2) \subset \bZ ^2$ 
is the triangle with vertices $(2,0)$, $(4,2)$, $(4,0)$. Indeed,
$\overline{G\cdot x_2}$ contains both $X_{12}$ and 
$Y_2 :=\overline{G\cdot [x^4 x' y']}$. Thus, the segments 
$[(2,0),(4,2)]$ and $[(4,2),(4,0)]$ are faces of $Q_2$. Since this
convex polytope consists of dominant weights, this leaves no other
choice.

By the description of their weight groups and moment polytopes,
$(X_1,L_1)$ and $(X_2,L_2)$ contain both $(X_{12},L_{12})$. So we may
glue them to a multiplicity-free PSSV $(X,L)$ with moment polytope 
$Q:= Q_1 \cup Q_2$, the triangle with vertices $(0,0)$, $(4,0)$,
$(4,2)$. The corresponding partial subdivision $\cQ$ consists of 
the triangles $Q_1$, $Q_2$, the segments $[(0,0),(2,0)]$,
$[(2,0),(4,2)]$, $[(4,2),(4,0)]$, and the points $(2,0)$, $(4,2)$.

We claim that there exists no spherical variety over $\bP(V)$ with
moment polytope $Q$. In other words, $Q$ is not the moment polytope of
any irreducible multiplicity-free $G$-subvariety 
$X\subseteq \bP(V)$. 
Indeed, if $X$ contains both closed $G$-orbits, then it 
contains $X_{12}$; this contradicts the fact that $Q_{12}$ is not a
face of $Q$. So $X$ contains a unique closed orbit. Hence one of the 
projections $\bP(V) -\to \bP(V(2,0))$, $\bP(V) - \to \bP(V(4,2))$ 
restricts to a finite morphism $X \to \bP(V(2,0))$ or 
$X \to \bP(V(4,2))$. The former case is excluded, since
$\rk(X)=2$. Thus, we reduce to the case where  
$X \subseteq \bP(V(4,2))$. Now $X$ is multiplicity-free of dimension
$4$ and contains semi-stable points, as $Q$ contains the origin. By
inspection, one obtains the unique possibility 
$X= \overline{G\cdot [x^2 y^2 x' y']}$. But then the moment polytope
contains $(0,2)$, a contradiction.

One checks that $\Aut^G(X_{12},L_{12})= \bG_m \times \bG_m$ and
likewise for $X_1$, $X_2$. It follows that 
$H^0(\cQ,\Aut) = \bG_m \times \bG_m$ and $H^1(\cQ,\Aut)$ is trivial. 
In other words, $\Aut^G(X,L)= \bG_m \times \bG_m$, and any PSSV
with building blocks $(X_1,L_1)$, $(X_2,L_2)$, $(X_{12},L_{12})$ is
isomorphic to $(X,L)$.

\medskip

\noindent
5) \emph{Simple PSSV's}. 

Recall that a $G$-variety is said to be \emph{simple} if it contains a
unique closed $G$-orbit. Let $(X,L)$ be a simple multiplicity-free
PSSV with closed orbit $Z$. Then the $G$-module $H^0(Z,L)$ is simple,
and hence lifts to a unique simple submodule of $H^0(X,L)$. Clearly,
the latter submodule is base-point-free and invariant under
$\Aut^G(X,L)$. Thus, there exist a unique dominant weight $\lambda$
and a finite $G$-morphism $f:X \to \bP(V(\lambda))$ such that
$L=f^*\cO(1)$. Moreover, $\Aut^G(X,L)$ is a finite extension of
$\bG_m$. Conversely, any PSSV over the projectivization of a simple
$G$-module is also simple, as follows, e.g., from Lemma \ref{cover}.

Given $\lambda\in\Lambda^+$, we claim that 
\emph{there are only finitely many multiplicity-free PSSV's $(X,f)$ over}
$\bP(V(\lambda))$. Indeed, 
there are only finitely many spherical varieties over
$\bP(V(\lambda))$ by Theorem \ref{finiteness}, so that we may fix the
building blocks $(Y,L_Y)$ of $(X,L)$. Then each $\Aut^G(Y,L_Y)$ is a
finite extension of $\bG_m$. Moreover, $H^1(\cQ,\bG_m)$ is trivial,
since the poset $\cQ$ admits a unique minimal element. It follows that
$H^1(\cQ,\Aut)$ is finite. By Lemma \ref{cohomology}, this implies our
claim.

\section{Families}

\subsection{Families of affine SSV's}

From now on we consider families of varieties over schemes; by a scheme, 
we mean a Noetherian scheme over our base field $k$ (still assumed to be
algebraically closed, of characteristic zero). Morphisms (resp. products)
are understood to be morphisms (resp. products) over $k$. 

Recall from \cite[Section 7]{AB-I} that a 
\emph{family of affine $G$-varieties} over a scheme $S$ is a morphism
$\pi: \cX \to S$ satisfying the following conditions :

\smallskip

\noindent
(i) $\cX$ is a scheme equipped with an action of the constant group
scheme $G\times S \to S$.

\smallskip

\noindent
(ii) $\pi$ is flat, affine, and of finite type, with geometrically
connected and reduced fibers.

\smallskip

A \emph{$G$-morphism} of families over the same scheme $S$ is an
equivariant morphism over $S$. The family $\pi:\cX \to S$ is
\emph{trivial} if it is $G$-isomorphic to $X\times S$ equipped with
the projection to $S$, where $X$ is some affine $G$-variety.

Given a family of affine $G$-varieties $\pi:\cX \to S$ and a point
$s\in S$, the geometric fiber $\cX_{\bar{s}}$ is an affine
$G(k(\bar{s}))$-variety, where $k(\bar{s})$ denotes an
algebraic closure of the residue field $k(s)$. Moreover, 
$$
\cR := \pi_*(\cO_{\cX})
$$ 
is a flat $\cO_S$-algebra equipped with a compatible action of $G$,
which is rational by \cite[I.1]{MFK}; we say that $\cR$ is a
$\cO_S$-$G$-\emph{algebra}. This yields an isomorphism of 
$\cO_S$-$G$-modules
\begin{equation}\label{OSG}
\cR \simeq \bigoplus_{\lambda \in \Lambda^+} 
\cF_{\lambda}\otimes V(\lambda),
\end{equation}
where each $\cF_{\lambda}$ is a flat $\cO_S$-module. Moreover, each
$\cF_{\lambda}$ is a coherent sheaf of modules over the invariant
subalgebra, $\cR^G \simeq \cF_0$. 

Thus, if $\cR^G$ is finitely generated over $\cO_S$ (e.g., if each
geometric fiber contains only finitely many orbits), then each
$\cF_{\lambda}$ is a locally free $\cO_S$-module, since we assume $S$
to be Noetherian. The rank of $\cF_{\lambda}$ at any $s\in S$ is the
multiplicity of the simple $G(k(\bar{s}))$-module 
$k(\bar{s}) \otimes V(\lambda)$ in the affine ring of~$\cX_{\bar{s}}$.

Taking $U$-invariants in (\ref{OSG}) yields an isomorphism of
$\cO_S$-$T$-algebras
\begin{equation}\label{OST}
\cR^U \simeq \bigoplus_{\lambda\in\Lambda^+} \cF_{\lambda}.
\end{equation}

\begin{definition}
Given a subgroup $\Gamma$ of $\Lambda$, a 
\emph{family of affine SSV's with weights in $\Gamma$} 
is a family of affine $G$-varieties $\pi: \cX \to S$ such that each
geometric fiber is a stable spherical variety with weights in $\Gamma$. 
\end{definition}

If $S$ is connected and some geometric fiber $X = \cX_{\bar{s}}$ is
multiplicity-free (e.g., a spherical variety over $k(\bar{s})$), then
the $\cO_S$-module $\cF_{\lambda}$ is invertible whenever
$\lambda \in \Lambda^+(X)$, and $\cF_{\lambda}=0$ otherwise. Thus, all
the geometric fibers are multiplicity-free SSV's with the same weight
set. 

We now obtain an important isotriviality result for families of
spherical varieties.

\begin{proposition}\label{isotrivial}
Let $\pi : \cX \to S$ be a family of affine spherical varieties over 
an excellent integral scheme. Then there exist a non-empty open
subscheme $S_0$ of $S$ and a finite surjective morphism 
$\varphi : S' \to S_0$ such that the pull-back family 
$\pi': \cX \times_S S'\to S'$ is trivial. 
\end{proposition}

\begin{proof}
By \cite[Lemma 7.4]{AB-I}, the family 
$$
\pi//U : \cX//U := \Spec_{\cO_S}\cR^U \to S
$$ 
is locally trivial, with fiber an affine toric variety $Y$. Replacing
$S$ with an open subset, we may assume that there is an isomorphism 
$\cX//U \to Y \times S$ over $S$. This yields a morphism 
$\varphi: S \to M_Y$, where $M_Y$ denotes the moduli scheme of affine
spherical varieties of type $Y$ defined in \cite[Section 1.3]{AB-III}; 
then $\cX$ is the pull-back of the universal family of $M_Y$. Recall 
from \cite[Corollary 3.4]{AB-III} that $M_Y$ has an action of $T$ with
finitely many orbits; these are in bijection with the isomorphism
classes of affine spherical varieties $X$ such that $X//U \simeq Y$. 
So, by shrinking $S$ again, we may assume that the image of $\varphi$
is contained in a unique $T$-orbit $\cO$. After a finite surjective
base change, we may assume that $\varphi: S \to \cO$ lifts to a
morphism $\psi: S \to T$. Since the $T$-action on $M_Y$ lifts to an
action on the universal family by \cite[Section 2.1]{AB-III}, it
follows that $\pi$ is trivial.
\end{proof}

Next we study \emph{one-parameter degenerations} of affine spherical
varieties, that is, families of affine $G$-varieties $\pi:\cX \to S$,
where $S$ is a regular integral scheme of dimension $1$, and the
geometric generic fiber $\cX_{\bar{\eta}}$ is a spherical variety. 

Here is a construction of such families. Let $\tX$ be an affine
spherical $\tG$-variety (recall that $\tG = \bG_m \times G$) equipped
with a surjective $\tG$-morphism $f:\tX \to \bA^1$, where $\tG$ acts on
$\bA^1$ via the trivial action of $G$ and the scalar action of
$\bG_m$. Then $f$ restricts to a trivial family over 
$\bA^1 \setminus \{0\} \simeq \bG_m$, with fiber some affine spherical
$G$-variety $X$. 

\begin{lemma}\label{standard}
With the preceding notation, let 
$\tC\subseteq \bR \times \Lambda_{\bR}$ (resp.
$C\subseteq \Lambda_{\bR}$) be the weight cone of $\tX$
(resp. $X$). 

\smallskip

\noindent
(i) There exists a unique function $h : C \to \bR$ such that 
\begin{equation}\label{graph}
\tC = \{(t,\lambda)\in \bR \times C ~\vert~ h(\lambda) \le t\}.
\end{equation}
Moreover, $h$ is lower convex, piecewise linear, and takes rational
values at all rational points of $C$.

\smallskip

\noindent
(ii) The special fiber $f^{-1}(0)$ is reduced if and only if 
$h$ takes integral values at all points of 
$\Lambda^+(X) = C \cap \Lambda(X)$. 
Then $f^{-1}(0)$ is a stable spherical variety, so that $f$ is a
one-parameter degeneration.

\smallskip

\noindent
(iii) There exists a positive integer $N$ such that the pull-back
family under the morphism $\bA^1 \to \bA^1$, $z' \mapsto z^{'N}=z$,
has reduced fibers.
\end{lemma}

\begin{proof}
(i) Let $R$ (resp. $\tR$) be the affine ring of $X$ (resp. $\tX$). 
Then $\tR \subseteq R[z,z^{-1}]$, as 
$\tX \setminus f^{-1}(0)\simeq X \times \bG_m$. So
$\tC \subseteq \bR \times C$. Moreover, $(1,0) \in \tC$ since 
$f\in \tR$, but $(-1,0) \notin \tC$ since $f$ is not invertible in
$\tR$. Since $\tC$ is a closed convex cone, it may be written as 
(\ref{graph}) for a unique convex function $h$. The piecewise
linearity and rationality of $h$ follow from the fact that $\tC$ is a
rational polyhedral cone.  

(ii) Assume that there exists $\lambda\in\Lambda^+(X)$ such that 
$h(\lambda) \notin \bZ$. Let $n$ be the smallest integer such that
$n > h(\lambda)$, and let $f_{\lambda}\in R$ be a $B$-eigenvector of
weight $\lambda$. Then $\tf_{\lambda} := z^n f_{\lambda}$ is in
$\tR$, but not in $z\tR$. Now let $N$ be a positive integer such that
$N h(\lambda) \in \bZ$. Then $\tf_{\lambda}^N := z^{Nn} f_{\lambda}^N$
is in $z\tR$, since $Nn \ge N h(\lambda) + 1$. So the algebra
$\tR/z\tR$ is not reduced.

Conversely, if $h$ takes integral values at all points of
$\Lambda^+(X)$, then the preceding argument shows that 
no $B$-eigenvector in $\tR/z\tR$ is nilpotent. It follows that 
$\tR/z\tR$ is reduced. On the other hand, $\tR/z\tR$ is
multiplicity-free. By Lemma \ref{multiplicityfree}, this implies
that the special fiber is a SSV.

(iii) The pull-back under consideration is a similar family, but where
$h$ is replaced with $Nh$. Since $h$ is rational, and linear on each
cone of a finite subdivision of $C$ into rational polyhedral convex
cones, we may choose $N$ so that $Nh$ takes integral values at all
points of $\Lambda^+(X)$. Then the assertion follows from (ii).
\end{proof}

A one-parameter degeneration as in Lemma \ref{standard} will be called
\emph{standard}. We now show that every one-parameter degeneration
becomes standard after appropriate base changes.

\begin{proposition}\label{model}
Let $\pi: \cX \to S = \Spec A$ be a one-parameter degeneration of
affine spherical varieties, where $A$ is a discrete valuation
ring. Then, after a finite surjective base change, $\pi$ is isomorphic
to the pull-back of a standard family $p:\tX \to \bA^1$.
\end{proposition}
 
\begin{proof}
We adapt the argument of \cite[Proposition 7.13]{AB-I}.
Let $\eta = \Spec K$ be the generic point of $S$, where $K$ is the
fraction field of $A$, and let $s = \Spec(A/zA)$ be the closed point,
where $z$ is a generator of the maximal ideal of $A$. 
Applying Proposition \ref{isotrivial}, we may assume (after a finite
surjective base change) that $\cX_{\eta} = X \times \{\eta\}$, where
$X = \Spec R$ is an affine spherical variety. Then 
$\cX_{\eta} = \Spec \cR_{\eta}$, where  
$$
\cR_{\eta} = K \otimes R \simeq  
\bigoplus_{\lambda\in \Lambda^+(X)} K \otimes V(\lambda).
$$
Moreover, $\cX = \Spec \cR$, where 
$\cR = \bigoplus_{\lambda\in \Lambda^+(X)} 
\cF_{\lambda} \otimes V(\lambda)$,
and each $\cF_{\lambda}$ is an invertible $R$-submodule of $K$. Thus, 
$$
\cR = \bigoplus_{\lambda\in \Lambda^+(X)} 
z^{h(\lambda)} R \otimes V(\lambda)
$$ 
for some function $h:\Lambda^+(X) \to \bZ$. 

We claim that the subspace
$$
\tR := \bigoplus_{\lambda\in \Lambda^+(X)} 
z^{h(\lambda)} k[z] \otimes V(\lambda) \subseteq \cR
$$
is in fact a subalgebra. Indeed, consider three weights $\lambda$,
$\mu$, $\nu$ in $\Lambda^+(X)$, such that $V(\nu)$ occurs in the
decomposition of the product $V(\lambda) \cdot V(\mu)$ in $R$. Then
the product $z^{h(\lambda)}V(\lambda)\cdot z^{h(\mu)} V(\mu)$ (in
$\cR$) contains $z^{h(\lambda) + h(\mu)} V(\nu)$, so that 
$h(\nu) \le h(\lambda) + h(\mu)$. This implies our claim.

Clearly, $\tR$ contains the polynomial ring $k[z]$, and 
$$
k(s) \otimes \tR/z\tR \simeq \cR/z\cR,
$$ 
where $k(s) = A/zA$. Since the algebra $\cR/z\cR$ is the affine ring
of the special fiber, it is finitely generated over $k(s)$. It follows
that $\tR/z\tR$, and hence $\tR$, is finitely generated as
well. Moreover, $\tR/z\tR$ is reduced, since $\cR/z\cR$ is. By
considering powers of $B$-eigenvectors as in the proof of Lemma
\ref{standard}, this implies the equality $h(n\lambda) = n h(\lambda)$
for any positive integer $n$ and any $\lambda \in \Lambda^+(X)$. So
$\tX := \Spec \tR$ is an irreducible multiplicity-free $\tG$-variety
with a saturated weight set, i.e., a spherical
$\tG$-variety. Moreover, the morphism $z:\tX \to \bA^1$ is a standard
degeneration. Finally, the multiplication map yields an isomorphism of
$A$-algebras  
$$
A \otimes_{k[z]} \tR \simeq \cR,
$$ 
which completes the proof.
\end{proof}

\subsection{Families of polarized SSV's}

A \emph{family of polarized $G$-varieties over $S$} is a
pair $(\pi:\cX \to S,\cL)$ satisfying the following conditions.

\smallskip

\noindent
(i) $\cX$ is a scheme equipped with an action of the constant group 
scheme $G\times S$ over $S$.

\smallskip

\noindent
(ii) $\pi$ is flat, proper, with geometrically connected and reduced
fibers. 

\smallskip

\noindent
(iii) $\cL$ is a $\pi$-ample, $G$-linearized invertible sheaf on
$\cX$.

\smallskip

We then put 
$$
\cR(\cX,\cL):= \bigoplus_{n=0}^{\infty} \pi_*(\cL^n).
$$
This is a sheaf of $\cO_S$-$\tG$-algebras, the 
\emph{(relative) section ring of $(\cX,\cL)$}.

\begin{definition}
Given a subgroup $\Gamma$ of $\tLambda$, a 
\emph{family of PSSV's with weights in $\Gamma$}
is a family of polarized $G$-varieties such that every geometric fiber
is a PSSV with weights in $\Gamma$.
\end{definition}

From now on we will only consider families of \emph{multiplicity-free}
PSSV's $\pi:\cX \to S$, with weights in a prescribed subgroup
$\Gamma$ of $\tLambda$. By Proposition \ref{pssv}, 
$H^1(\cX_{\bar{s}},\cL^n_{\bar{s}})=0$ 
for all geometric points $\bar{s}$ and for all $n\ge 1$. Thus, every
$\cO_S$-module $\pi_*(\cL^n)$ is locally free and satisfies
$\pi_*(\cL^n) \otimes_{\cO_S} k(\bar{s}) = 
H^0(\cX_{\bar{s}},\cL^n_{\bar{s}})$
for all $\bar{s}$, by the theorem on cohomology and
base change \cite[Chapter III, Theorem 12.11]{Ha}. This yields 
isomorphisms of $\cO_S$-$G$-modules
\begin{equation}\label{dir}
\pi_*(\cL^n) \simeq \bigoplus_{\lambda \in \Lambda^+} 
\cF_{n,\lambda}\otimes V(\lambda),
\end{equation}
where each non-zero $\cF_{n,\lambda}$ is an invertible $\cO_S$-module. 

In other words, the morphism
$$
\tpi: \tcX := \Spec_{\cO_S} \, \cR(\cX,\cL) \to S
$$
is a family of affine multiplicity-free SSV's (for $\tG$)
with $\cX = \Proj_{\cO_S} \cR(\cX,\cL)$ and 
$\cL = \cO_{\cX}(1)$. Moreover, 
$$
\cR(\cX,\cL) \simeq \bigoplus_{\tlambda \in \tLambda^+}
\cF_{\tlambda} \otimes V(\tlambda), \quad
\cR(\cX,\cL)^U \simeq \bigoplus_{\tlambda \in \tLambda} \cF_{\tlambda},
$$
and $\cR(\cX,\cL) \otimes_{\cO_S} k(\bar{s}) = 
R(\cX_{\bar{s}},\cL_{\bar{s}})$
for any geometric point $\bar{s}$. The sheaves 
$\cF_{\tlambda}$ are equipped with multiplication maps
\begin{equation}\label{mult}
m_{\tlambda,\tmu} : \cF_{\tlambda}\otimes_{\cO_S} \cF_{\tmu} 
\to \cF_{\tlambda + \tmu}~.
\end{equation}
In particular, this yields $N$-th power maps
\begin{equation}\label{power}
\cF_{\tlambda}^N \to \cF_{N\tlambda}
\end{equation}
which are isomorphisms whenever $\cF_{\tlambda}$ is non-zero,
since the fibers of $\pi$ are reduced.

Thus, if $S$ is connected then all the geometric fibers have the same
moment set $Q$; we say that the family is \emph{of type $Q$}. 

Next we study the families of SSV's over the projectivization of a
fixed $G$-module $V$.

\begin{definition}
A \emph{family of SSV's over $\bP(V)$} is a pair 
$(\pi:\cX \to S, f:\cX \to \bP(V))$, where 
$(\cX,\cL:=f^*\cO(1))$ is a family of PSSV's. 

Then the product morphism $f \times \pi: \cX \to \bP(V) \times S$ 
is finite, so that every geometric fiber of $\pi$ is a SSV over
$\bP(V)$. We also say that $\cX$ is a \emph{SSV over} 
$\bP(V)\times S$.

A \emph{morphism} between two SSV's over $\bP(V)\times S$ is a 
$G$-morphism of schemes over $\bP(V)\times S$.
\end{definition}

We now obtain two boundedness results for these families, which will
play an essential role in the construction of the moduli space.

\begin{lemma}\label{multiplication}
Given a finite-dimensional $G$-module $V$, there exists a finite
collection $\cQ$ of rational convex polytopes in $\Lambda^+_{\bR}$
such that the following properties hold for any SSV $\cX$ over
$\bP(V)\times S$: 

\smallskip

\noindent
(i) The moment polytope of any $G$-subvariety of any geometric fiber 
$\cX_{\bar{s}}$ is a union of polytopes in $\cQ$. 

\smallskip

\noindent
(ii) For any $Q \in \cQ$ and $\tlambda$, $\tmu$ in 
$\Gamma \cap \Cone(Q)$, the map $m_{\tlambda,\tmu}$ 
of (\ref{mult}) is an isomorphism.

\end{lemma}

\begin{proof}
By Theorem \ref{finiteness}, the set of moment polytopes of 
spherical varieties over $\bP(V)$ is finite. Therefore, we may
choose a finite common subdivision $\cQ$ of all these polytopes
by rational convex polytopes. Clearly, this subdivision satisfies (i).

Let $s\in S$ and put $X := \cX_{\bar{s}}$, $L := \cL_{\bar{s}}$. Then,
by multiplicity-freeness, the multiplication map
$R(X,L)^U_{\tlambda} \otimes R(X,L)^U_{\tmu} \to 
R(X,L)^U_{\tlambda+\tmu}$
is an isomorphism whenever both $\tlambda$, $\tmu$ belong to 
$\Gamma \cap \Cone Q(Y,L)$, where $Y$ is an irreducible component of
$X$. It follows that $\cQ$ satisfies (ii).
\end{proof}

Any family of SSV's $f\times \pi :\cX \to \bP(V)\times S$ 
defines a morphism of $\cO_S$-$\tG$-algebras
$$
f^* : \cO_S \otimes \Sym(V^*) \to \cR(\cX,f^*\cO(1)) := \cR,
$$
which makes $\cR$ a finite module over $\cO_S\otimes \Sym(V^*)$. 
It follows that $\cR^U$ is a finite module over 
$\cO_S\otimes \Sym(V^*)^U$. For every $\tlambda\in\tLambda^+$, we
denote by 
\begin{equation}\label{pullback}
f^*_{\tlambda}: \cO_S\otimes \Sym(V^*)^U_{\tlambda} 
\to  \cR^U_{\tlambda} = \cF_{\tlambda}
\end{equation}
the restriction of $f^*$ to the component of weight $\tlambda$.

\begin{lemma}\label{generation}
Let $\cQ$ be as in Lemma \ref{multiplication}. Then:

\smallskip

\noindent
(i) There exists a finite subset $F$ of $\tLambda^+$ such that the
monoid $\Gamma \cap \Cone(Q)$ is generated by a subset of $F$, for
any $Q \in \cQ$.  

\smallskip

\noindent
(ii) For any SSV $\cX$ over $\bP(V)\times S$, the $\cO_S$-algebra  
$\cR^U$ is generated by the $\cF_{\tlambda}$, where $\tlambda \in F$.
Moreover, there exists a positive integer $N$ (depending only
on $V$) such that the map $f^*_{N\tlambda}$ of (\ref{pullback})
is surjective for any $\tlambda \in \tLambda^+$.
\end{lemma}

\begin{proof}
By Gordan's lemma, every monoid $\Gamma \cap \Cone(Q)$, $Q\in \cQ$,
is finitely generated. Thus, we may choose a finite subset $F\subseteq
\tLambda^+$ containing generators of all these monoids. Then the
sheaves $\cF_{\tlambda}$, $\tlambda \in F$, generate the algebra
$\cR^U$ by Lemma \ref{multiplication} (ii). 

To complete the proof, by Lemma \ref{multiplication} again, it
suffices to show the existence for any $\tlambda$ of $N = N(\tlambda)$ 
such that $f^*_{N\tlambda}$ is surjective. Let
$$
\cR^U_{(\tlambda)} := \bigoplus_{n=0}^{\infty}
\cR(\cX,\cL)^U_{n\tlambda}
=\bigoplus_{n=0}^{\infty} \cF_{n\tlambda}.
$$
This is a finitely generated, graded $\cO_S$-subalgebra of $\cR^U$; we
may assume that $\cR^U_{(\tlambda)} \neq \cO_S$. 
Let $N_0$ be the smallest positive integer 
such that $\cF_{N_0\tlambda}\ne 0$. By saturation,  
$\cF_{n\tlambda}\ne 0$ if and only if $N_0$ divides $n$. 
It follows that 
$$
\cR^U_{(N\tlambda)} \simeq \Sym_{\cO_S}(\cF_{N\tlambda})
$$
whenever $N$ is a multiple of $N_0$.
In particular, $\cR^U_{(N\tlambda)}$ is locally isomorphic to
the polynomial ring $\cO_S[z]$, where $z$ is a variable of
degree $N$.

Define likewise the finitely generated, graded algebra 
$$
\Sym(V^*)^U_{(\tlambda)} = \bigoplus_{n=0}^{\infty} 
\Sym(V^*)^U_{n\tlambda}.
$$
Then there exists a positive integer $N_1$ such that the algebra 
$\Sym(V^*)^U_{(N\tlambda)}$ is generated by its subspace of 
degree $N$, whenever $N$ is a multiple of $N_1$. On the other hand, by
weight considerations, $\cR^U_{(N\tlambda)}$ is a finite module over  
$\cO_S \otimes \Sym(V^*)^U_{(N\tlambda)}$ for any $N$.
It follows easily that the desired surjectivity holds for any $N$
divisible by $N_0$ and $N_1$.
\end{proof}

\begin{remark}
By similar arguments, one obtains the existence of a positive integer
$N'$ (depending only on $V$) such that the map
$$
f^*_n : \cO_S \otimes \Sym^n(V^*) \to \cR_n
$$ 
is surjective for any SSV $f\times \pi: \cX \to \bP(V) \times S$ and
for any multiple $n$ of $N'$. In particular, $f^*\cO(N')$ is very
ample relatively to $\pi$. Since $f^*\cO(1)$ is $\pi$-globally
generated, it follows that $f^*\cO(n)$ is $\pi$-very ample for any
$n\ge N'$. 
\end{remark}

\subsection{Families of stable spherical pairs}
\label{sec:Families of stable spherical pairs}

In this section, we introduce the notion of stable spherical pairs,
and we reduce their classification to that of stable spherical
varieties over the projective spaces of certain $G$-modules.

We still consider families of \emph{multiplicity-free}
PSSV's with weights in a prescribed group $\Gamma$.

\begin{definition}
A \emph{family of stable spherical pairs} (in short, SSP's)
over a scheme $S$ consists of a family ($\pi:\cX \to S, \cL)$ 
of PSSV's, together with a section $\sigma \in H^0(\cX,\cL)$ 
satisfying the following condition:

\smallskip

\noindent
For any geometric point $\bar{s}$, the pull-back of $\sigma$
to the geometric fiber $\cX_{\bar{s}}$ does not vanish on any 
orbit of $G(k(\bar{s}))$.

\smallskip

Then the divisor of zeroes $\cD := (\sigma)_0$ 
is an effective, $\pi$-ample Cartier divisor on $\cX$. Since $\cL =
\cO_{\cX}(\cD)$, the pair $(\cX,\cD)$ encodes our triple 
$(\pi: \cX \to S,\cL,\sigma)$.

\smallskip

A \emph{morphism} from another pair $(\pi': \cX' \to S,\cD')$ to 
$(\pi: \cX \to S, \cD)$ is a $G$-morphism $\varphi: \cX' \to \cX$ over
$S$ such that $\cD' = \varphi^*(\cD)$.
\end{definition}

We now associate with any family of SSP's $(\pi:\cX \to S,\cD)$ 
of a fixed type $Q$, a family of SSV's over the
projectivization of a $G$-module $V_Q$ depending only on $Q$. 
By (\ref{dir}), we have a canonical isomorphism of $\cO_S$-$G$-modules
$$
\pi_*(\cL) = \bigoplus_{\lambda \in F} 
\cF_{1,\lambda} \otimes V(\lambda),
$$
where $F := \Gamma \cap Q$ is a finite set of dominant weights. 
Then 
$$
\sigma\in H^0(\cX,\cL) = \bigoplus_{\lambda \in F} 
H^0(S,\cF_{(1,\lambda)}) \otimes V(\lambda).
$$
Write accordingly
$\sigma = \sum_{\lambda\in F} f_{\lambda} \otimes x_{\lambda}$. 
This yields a linear map
$$
\gamma: \bigoplus_{\lambda\in F} \End  V(\lambda) \to 
H^0(\cX,\cL), \quad 
\sum u_{\lambda}\mapsto \sum 
f_{\lambda} \otimes u_{\lambda}(x_{\lambda}).
$$
In particular, $\gamma(\sum_{\lambda\in F} \id_{\lambda}) = \sigma$ 
with obvious notation. Also, $\gamma$ is $G$-equivariant, where $G$
acts on each 
$\End V(\lambda) \simeq V(\lambda)^* \otimes V(\lambda)$ 
via its action on $V(\lambda)$. Put
\begin{equation}\label{end}
V_Q := \bigoplus_{\lambda\in F} (\End V(\lambda))^*
= \bigoplus_{\lambda\in F} V(\lambda)\otimes V(\lambda)^*,
\end{equation}
regarded as a $G$-module via the action on the spaces $V(\lambda)^*$.
Then we have a morphism of $G$-modules
$$
\gamma:V^*_Q \to H^0(\cX,\cL).
$$

\begin{proposition}\label{correspondence}
Let $(\pi:\cX \to S,\cD)$ be a family of SSP's of type $Q$. Then
the above morphism $\gamma$ yields a $G$-morphism 
$f:\cX \to \bP(V_Q)$ such that $\cL = f^* \cO(1)$ and 
$\sigma = f^*(\sum_{\lambda\in F} \id_\lambda)$. In particular,
$\cX$ is an SSV over $\bP(V_Q)\times S$.

This defines a bijective correspondence from families of SSP's 
of type $Q$ over $S$, to SSV's of the same type over 
$\bP(V_Q)\times S$; this correspondence preserves morphisms.
\end{proposition}

\begin{proof}
Let $\bar{s}$ be a geometric point of $S$. Then, by assumption,
the translates $g \cD_{\bar{s}}$, $g\in G(k(\bar{s}))$, have no common
zero on the geometric fiber $\cX_{\bar{s}}$. Since $G(k)$ is dense in
$G(k(\bar{s}))$, the same holds for the translates 
$g \cD_{\bar{s}}$, $g\in G(k)$. Let $g_{\lambda}$ denote the image of $g$
in $\GL(V(\lambda)) \subset \End V(\lambda)$. Then
$$
\gamma(\sum_{\lambda \in F} g_{\lambda}) =  
\sum_{\lambda\in F} g_{\lambda} \sigma_{\lambda} = g\sigma.
$$
Thus, the subspace $\gamma(V^*_Q)$ of $H^0(\cX,\cL)$ 
is base-point-free in each geometric fiber. It follows that
$f:\cX \to \bP(V_Q)$ is well-defined, and $\cL = f^*\cO(1)$.
Clearly, $f^*(\sum_{\lambda\in F} \id_{\lambda}) = \sigma$.

Conversely, let $(\pi : \cX \to S, ~f : \cX \to \bP(V_Q))$ be 
a family of SSV's over $\bP(V_Q)$. Let 
$\sigma := f^*(\sum_{\lambda \in F} \id_{\lambda})$; this is 
a global section of $\cL := f^*\cO(1)$. We show that $\sigma$
does not vanish identically on any $G(k(\bar{s})$-orbit $Y$ in a
geometric fiber $X =\cX_{\bar{s}}$. We may assume that $Y$ is closed,
and (to simplify notation) $k(\bar{s}) = k$. Then 
$Y \simeq G/P$, where $P$ is a parabolic subgroup of $G$ 
containing $B$. By assumption, we have a finite morphism
$$
f_{\bar{s}} : X \to \bP(\bigoplus_{\lambda \in F}
V(\lambda) \otimes V(\lambda)^*).
$$
By highest weight theory, there exists a unique 
$\lambda = \lambda(Y) \in F$ such that $f_{\bar{s}}$ restricts to an 
embedding 
$$
Y \to \bP(V(\lambda) \otimes V(\lambda)^*),
\quad gP \mapsto [v \otimes gv_{\lambda^*}],
$$
where $v\in V(\lambda)$ is non-zero, and 
$v_{\lambda^*} \in V(\lambda)^*$ is a highest weight vector.
Moreover, 
$\sigma(v\otimes gv_{\lambda^*}) =
\langle v, gv_{\lambda^*} \rangle$
is non-zero for some $g\in G$, since the translates 
$gv_{\lambda^*}$ span $V(\lambda)^*$.

This establishes the desired correspondence, which is clearly
functorial. 
\end{proof}

For example, if $G = T$ is a torus, then the families of stable toric
pairs of type $Q$ may be identified with the families of stable toric
varieties over $\bP(V_Q)$, where
$V_Q := \bigoplus_{\lambda \in \Gamma \cap Q} k_{-\lambda}$. 
Here $k_{\lambda}$ denotes the line $k$ on which $G$ acts with
weight $\lambda$.

\begin{remark}
We mention the following connection with the singularities of pairs,
generalizing results of \cite[Section 5]{AB-II} which mainly
concern stable reductive varieties.

Recall that a spherical variety $X$ has two kinds of group boundaries:
$\partial_G X$, the codimension one part of the complement of the open
$G$-orbit, and $\partial_B X$, the complement of the open $B$-orbit
minus $\partial_G X$. A canonical divisor for $X$ is 
\begin{equation}\label{canonical}
K_X = - \Delta_G - \Delta_B,
\end{equation}
where $\Delta_G$ is the divisor $\partial_G X$ with reduced structure,
and $\Delta_B$ is a unique effective divisor with support 
$\partial_B X$. 

Next consider a multiplicity-free stable spherical variety $X$ with
convex moment set. Then $X$ is Cohen--Macaulay by Proposition
\ref{pssv}, and one easily shows that $X$ has only simple
crossings in codimension one. It follows that (\ref{canonical}) still
holds, where $\Delta_G$ denotes the sum of all  irreducible
$G$-invariant divisors which are not contained in the double locus,
and $\Delta_B$ denotes the sum of the $\Delta_B$'s of these components.

In \cite[Theorem 5.3]{AB-II} we proved that for a spherical
variety $X$ the pair $(X, \Delta_G + |\Delta_B|)$ has log
canonical singularities. (Here $|\Delta_B|$ means that one has
to pick a general element of this linear system.) And an easy
extension of \cite[Theorems 5.9, 5.12]{AB-II} gives the following:

\begin{proposition}
Let $(X,D)$ be a multiplicity-free SSP whose moment set $Q$ is
convex. Then for $0 < \varepsilon \ll 1$, the pair 
$(X, \Delta_G + |\Delta_B| + \varepsilon D)$ has semi-log canonical
singularities (resp. log canonical if $X$ is irreducible).
\end{proposition}

\end{remark}

\section{Moduli}

\subsection{Existence of a quasiprojective moduli scheme}

We fix a subgroup $\Gamma$ of $\tLambda$ and a finite-dimensional
$G$-module $V$. Consider the contravariant functor
$$
\cM = \cM_{\Gamma,\bP(V)}
$$ 
from schemes to sets, that assigns to any scheme $S$ the set of
isomorphism classes of multiplicity-free SSV's over $\bP(V) \times S$
with weights in $\Gamma$. 

\begin{theorem}\label{coarse}
The functor $\cM$ is coarsely represented by a quasiprojective 
scheme $M_{\Gamma,\bP(V)}$.
\end{theorem}
 
\begin{proof}
Choose a finite collection $\cQ$ of polytopes, a finite subset $F$ of
$\tLambda^+$, and a positive integer $N$ satisfying the statements of
Lemma \ref{generation}. 

For simplicity, we begin with the case where $N=1$. Then, for any
family
$f \times \pi : \cX \to \bP(V) \times S$ with section ring 
$\cR := \cR(\cX,f^*\cO(1))$, the maps $f^*_{\tlambda}$ of
(\ref{pullback}) are all surjective. Thus, the map
$$ 
f^*:\cO_S \otimes \Sym(V^*) \to \cR
$$
is surjective as well. In other words, the morphism
$$
\tf \times \tpi : \tcX \to \bA(V) \times S
$$
is a closed immersion, where $\tcX = \Spec_{\cO_S} \, \cR$.
The image of this morphism yields a family of $\tG$-subvarieties of
$\bA(V)$, parametrized by $S$. The corresponding Hilbert function 
\begin{equation}\label{hilbert}
h : \tLambda^+ \to \bN
\end{equation}
(in the sense of \cite[Definition 1.4]{AB-III})
is given by: $h(\tlambda) = 1$ if there exists $Q \in \cQ$ such that
$\tlambda \in \Gamma \cap \Cone(Q)$, and $h(\tlambda) = 0$
otherwise.

Recall from \cite{AB-III} that the families of $\tG$-subschemes of
$\bA(V)$ with Hilbert function $h$ admit a fine moduli space: the
invariant Hilbert scheme 
$$
H:=\Hilb_h^{\tG}(V),
$$
a quasi-projective scheme. Moreover, the families of
$\tG$-subvarieties of $\bA(V)$ with Hilbert function $h$ are
parametrized by an open subscheme $H'$ of $H$, the locus where the
fibers of the universal family are reduced; see
\cite[p.~265]{AB-II}. So we obtain a morphism $S \to H'$ such
that $\tcX$ is isomorphic to the pull-back of the universal family.
Conversely, the universal family over $H'$ is a family of stable
spherical $\tG$-subvarieties of $\bA(V)$ (by Lemma
\ref{multiplicityfree}), and hence a family of SSV's over
$\bP(V)$. Thus, $\cM$ is represented by the quasi-projective
scheme $H'$.

We now turn to the general case, where $N$ is arbitrary. 
A little problem we need to overcome is that some of our
generating sections of $\cR$ provided by Lemma \ref{generation} have
weights $\tlambda_i=(n_i,\lambda_i)$ of degrees $n_i>1$, and moreover 
only their $N$-th power may belong to the image of
$\cO_S\otimes\Sym(V^*)$. We solve this by taking $N$-th roots, thus
making extra choices, and then dividing by these choices. As a result,
the moduli space is only coarse and not fine.

The $N$-th roots live on a finite abelian Galois covering of
the base scheme $S$.  
Specifically, write $F = \{ \tlambda_1,\dots,\tlambda_m\}$ and
$\tlambda_i = (n_i,\lambda_i)$ for $i=1,\ldots,m$, where $n_i$
is a positive integer and $\lambda_i$ is a dominant weight.  Choose a
basis $(\varphi_{ij})_{j \in J_i}$ of the finite-dimensional vector
space
$$
\Sym(V^*)^U_{N\tlambda_i} = \Sym^{Nn_i}(V^*)^U_{N\lambda_i}.
$$
Consider again a family 
$f \times \pi : \cX \to \bP(V) \times S$
with section ring $\cR$, and the associated maps $f^*_{\tlambda}$ of
(\ref{pullback}). By Lemma \ref{generation}, each map
$f^*_{N\tlambda_i}$ is surjective. Let
$$
\sigma_{ij} := f^*_{N\tlambda_i}(\varphi_{ij}) \in 
H^0(S,\cF_{N\tlambda_i}),
$$
then each invertible $\cO_S$-module $\cF_{N\tlambda_i}$ is
generated by its global sections $\sigma_{ij}$, $j\in J_i$. 

Let $S'$ be the scheme obtained from $S$ by taking the $N$-th roots of
these sections for all $i,j$ (see \cite[\S 3]{EV}) and let $A$ be the 
corresponding product of the groups $\bmu_N$ of $N$-th roots of
unity. Then the finite abelian group $A$ acts on $S'$ with a flat
quotient map $p:S' \to S$. Moreover, every invertible
$\cO_{S'}$-module 
$$
\cF'_{\tlambda_i} := p^*\cF_{\tlambda_i}
$$ 
is equipped with sections $\sigma'_{ij}$ such that 
$(\sigma'_{ij})^N = \sigma_{ij}$. Since
$(\cF'_{\tlambda_i})^N$ is identified with
$\cF'_{N\tlambda_i}$ via
(\ref{power}), it follows that the $\cO_{S'}$-module
$\cF'_{\tlambda_i}$ is generated by the $\sigma'_{ij}$, $j\in
J_i$.
So, by Lemma \ref{generation} again, the $\cO_{S'}$-algebra
$$
\cR' := \cR \otimes_{\cO_S}\cO_{S'} \simeq 
\bigoplus_{\tlambda \in \tLambda^+} 
\cF'_{\tlambda} \otimes V(\tlambda)
$$
is generated by its subspaces 
$\sigma'_{ij} \otimes V(\tlambda_i)$. This defines a
$\tG$-module map
$$
\gamma : V \oplus \bigoplus_{i,j\in J_i} V(\tlambda_i) \to \cR',
$$
where the restriction $V \to \cR'$ is the composed map
$V \to H^0(S,\cR) \to \cR'$. The image of $\gamma$ generates the
$\cO_{S'}$-algebra $\cR'$. Let 
$$
V' := V^* \oplus \bigoplus_{i,j} V(\tlambda_i)^*
\quad \text{and}\quad
\tcX' := \Spec_{\cO_{S'}} \cR'.
$$ 
Then $\gamma$ yields a closed $\tG$-embedding 
$$
\tcX' \hookrightarrow \bA(V') \times S
$$ 
over $S'$. Moreover, the group $A$ acts linearly on $\bA(V')$ via the
trivial action on $V$ and the scalar action of each factor $\bmu_N$ on
the  corresponding factor $V(\tlambda_i)$; this action commutes with 
the $\tG$-action. The structure map $\tcX' \to S'$ is
$A$-equivariant, the projection $\tcX' \to \bA(V) \times S'$ is
finite, and the quotient 
$\tcX'/A$ is $\tG$-isomorphic to $\tcX$ (since 
$(\cR')^A = (\cO_{S'}\otimes_{\cO_S} \cR)^A \simeq \cR$).

Let $h:\tLambda^+ \to \bN$ be as in (\ref{hilbert}). Consider the
invariant Hilbert scheme $H := \Hilb_h^{\tG}(V')$, and its locus $H'$
parametrizing those families $\cX' = \Spec_{\cO_{S'}}(\cR')$ having
reduced fibers and such that the maps
$$
f^{'*}_{N\tlambda_i}: \cO_{S'} \otimes \Sym(V^*)^U_{N\tlambda_i} 
\to \cR^{'U}_{N\tlambda_i}
$$
are surjective for $i=1,\ldots,m$. Clearly, $H'$ is an open
subscheme, invariant under the $A$-action on $H$ (induced by the
linear $A$-action on $\bA(V')$ that commutes with the
$\tG$-action). Since $H$ is quasi-projective, the quotient $H'/A$ is
quasi-projective as well. 

We now show that $H'/A$ coarsely represents $\cM$. Indeed, we just saw
that every family of SSV's over $\bP(V)\times S$ is the pullback of
the universal family over $H'$ under a functorially defined 
$A$-equivariant morphism $S' \to H'$, $S=S'/A$.  This gives a
morphism $\phi:\cM \to \Hom(*,H'/A)$ to the categorical quotient
$H'/A$, and $\phi$ is universal among such morphisms.

The universal family over $H'$ is a family of affine SSV's for
$\tG$ (by Lemma \ref{multiplicityfree}) which is finite over 
$\bA(V) \times H'$ (as follows from Lemma \ref{generation}). Thus,
each closed point of $H'/A$ defines a SSV over $\bP(V)$, which is
unique up to $G$-isomorphism.  Hence, $H'/A$ satisfies the conditions
in the definition of the coarse moduli space of $\cM$.
\end{proof}

\begin{remarks}\label{remarks}
1) Consider the subfunctors $\cM_{\Gamma,\bP(V),Q}$ of
$\cM_{\Gamma,\bP(V)}$, obtained by prescribing the moment set $Q$. 
Clearly, each $\cM_{\Gamma,\bP(V),Q}$ is an open and closed
subfunctor of $\cM_{\Gamma,\bP(V)}$. This yields moduli schemes
$M_{\Gamma,\bP(V),Q}$ which are unions of connected components of 
$M_{\Gamma,\bP(V)}$. 

\medskip

\noindent
2) Likewise, fix a subdivision $\cQ$ of $Q$ and consider the subfunctor
$\cM_{\Gamma,\bP(V),\cQ}$ where the geometric fibers have all moment set
$Q$ and subdivision $\cQ$. Then $\cM_{\Gamma,\bP(V),\cQ}$ is a locally
closed subfunctor. Indeed, it is given by the conditions that the map
$m_{\tlambda,\tmu}$ of (\ref{mult}) is an isomorphism whenever
$\tlambda$, $\tmu$ belong to the same cone over a polytope in $\cQ$, and
is zero otherwise. Further, it suffices to check these conditions for
finitely many pairs $(\tlambda,\tmu)$, by Lemmas \ref{multiplication}
and \ref{generation}. 

\medskip

\noindent
3) It is easy to check that any morphism of SSV's over 
$\bP(V)\times S$ of the same type $Q$ is an isomorphism, and
the automorphism group functor of any such SSV is representable by a
finite group scheme. Together with the argument of Theorem
\ref{coarse}, it follows that the families of SSV's over $\bP(V)$ with
weights in $\Gamma$ and type $Q$ form an algebraic stack of finite
type in the sense of \cite[Chapitre 4]{LM}, with coarse moduli space
$M_{\Gamma,\bP(V),Q}$.  In fact, this stack is proper, as will follow
from Proposition \ref{DVR}.
\end{remarks}

\subsection{Projectivity}

In this section, we show that the schemes $M_{\Gamma,\bP(V)}$
of Theorem \ref{coarse} are projective. Since we know that they are
quasiprojective, it suffices to check that the valuative
criterion of properness is satisfied. This is the content of the
following result. 

\begin{proposition}\label{DVR}
Let $A$ is a discrete valuation ring and let $S = \Spec A$ with
generic point $\eta$. Then every SSV $\cX_{\eta}$ over
$\bP(V) \times \eta$ extends to a unique SSV $\cX$ over 
$\bP(V) \times S$, possibly after a finite surjective base change 
$S' \to S$. 
\end{proposition}

\begin{proof}
Let $s=\Spec(A/zA)$ be the closed point of $S$, where $z$ is a
generator of the maximal ideal of $A$; let $K=k(\eta)$ be the fraction
field of $A$. Consider the section ring $\cR_{\eta}$ of $\cX_{\eta}$;
this is a $K$-$\tG$-algebra, finite over $K\otimes \Sym(V^*)$. The
desired extension $\cX$ corresponds to an $A$-$\tG$-subalgebra $\cR$
of $\cR_{\eta}$, finite over $A\otimes \Sym(V^*)$ and such that
$\cR_s = \cR/z\cR$ is reduced.

First we consider the case where the geometric generic fiber
$\cX_{\bar{\eta}}$ is a spherical variety. Then $\cR$ is an integral
domain having the same quotient field as $\cR_{\eta}$. Further, since
$\cR/z\cR$ satisfies ($R_0$) and ($S_1$), then $\cR$ satisfies ($R_1$)
and ($S_2$). Thus, $\cR$ is normal by Serre's criterion : so $\cR$
is the integral closure of $A[V]$ in $\cR_{\eta}$. This shows the
uniqueness of the extension $\cX$.

For the existence, we consider of course the integral closure $\cR$ of
$A[V]$ in $\cR_{\eta}$. This $A$-$\tG$-algebra is clearly finite over
$A[V]$ and flat over $A$, and hence multiplicity-free. It remains to
show that $\cR/z\cR$ (the section ring of the special fiber) is reduced,
possibly after a finite surjective base change. For this, we may
assume (by Proposition \ref{isotrivial}) that  
$$
\cR_{\eta} = K \otimes R \simeq 
\bigoplus_{\tlambda \in\tLambda(X)} K \otimes V(\lambda),
$$ 
where $R$ is the affine ring of an affine spherical $\tG$-variety
$X$. As in the proof of Proposition \ref{model}, it follows that 
$\cR$ is obtained by base change from a standard family 
$z:\tX \to \bA^1$. The special fiber of this family may be non-reduced,
but it becomes reduced after a finite base change $z = z^{'N}$, by
Lemma \ref{standard}.

We now show that our extension $\cX$ is compatible with
$G$-subvarieties. Let $\cY_{\bar{\eta}}$ be an irreducible
$G$-subvariety of $\cX_{\bar{\eta}}$. After a finite extension, we may
assume that $\cY_{\bar{\eta}}$ is obtained by base change from a
$G$-subvariety $\cY_{\eta}$ of $\cX_{\eta}$. Let $\cY$ be the closure
of $\cY_{\eta}$ in $\cX$. Then $\cY$ is flat over $S$ and
finite over $\bP(V)\times S$. Moreover, its special fiber $\cY_s$ is a 
closed $G$-subscheme of $\cX_s$, and the weight set of its affine cone
satisfies
$$
\tLambda^+(\tcY_s) = \tLambda^+(\tcY_{\eta}) = \tLambda^+(\tcY).
$$
Thus, $\tLambda^+(\tcY_s)$ is the intersection of $\tLambda^+(\tcX)$
with a face of the corresponding weight cone. By multiplicity-freeness
(as in the proof of Lemma \ref{reduced}), it follows that $\tcY_s$ is
reduced. So $\cY_s$ is reduced as well, and $\cY$ yields the desired
extension. 

Finally, we treat the general case, where $\cX_{\bar{\eta}}$ may have
several irreducible components. After a finite extension, we may
assume that all the irreducible $G$-subvarieties of $\cX_{\bar{\eta}}$
are defined over $K$. By Propositions \ref{ssv} and \ref{affinecone}, 
this yields affine spherical $\tG$-varieties $\tcY_{\eta}$ such that 
$\tcX_{\eta}= \varinjlim \tcY_{\eta}$. By the first step of the proof,
after a further finite base change, each $\tcY_{\eta}$ extends uniquely
to a family $\tcY$, finite over $V\times S$. Since these extensions
are compatible with $\tG$-subvarieties, they form a directed system of
SSV's. The direct limit of this system is the desired extension. 
\end{proof}

Next we consider the moduli scheme 
$$
M := M_{\Gamma,\bP(V),Q},
$$ 
where $Q$ is assumed to be convex. Then any SSV $(X,f)$ of type $Q$ is 
equidimensional, as follows, e.g., from Proposition \ref{pssv} (vii).
Clearly, the associated cycle $f_*[X]$ in $\bP(V)$ depends only on the 
isomorphism class of $(X,f)$, and its dimension and degree depend only
on the combinatorial data $\Gamma$, $Q$. This defines a cycle map
from the underlying set of $M$ to some Chow variety, $\Chow \bP(V)$
(see \cite[Chapter 1]{Ko} for background on Chow varieties and related
issues of semi-normality). 

We now show that this cycle map, properly defined, is a finite morphism.
(This generalizes a result of \cite[Section 4.4]{AB-II} regarding stable
reductive varieties, but the exposition there is inaccurate.)

\begin{proposition}\label{chow}
With the preceding notation, the cycle map extends to a finite
morphism $\gamma : M^{sn} \to \Chow \bP(V)$, where $M^{sn}$ denotes
the semi-normalization of $M$.
\end{proposition}

\begin{proof}
To show the existence of $\gamma$, we construct a kind of ``universal
family'' of cycles over $M^{sn}$. For this, we follow the 
notation as in the proof of Theorem \ref{coarse}. The subgroup 
$\bG_m = \bG_m\times \{1\}$ of $\tG = \bG_m \times G$ acts 
linearly on $\bA(V')$ with positive weights. Let $\bP(V')$ be the
associated weighted projective space, which may differ from the
projectivization of $V'$, but contains the usual projective
space $\bP(V)$. Then $\cX' = \Proj_{\cO_{S'}}\cR'$ is a
$G$-subscheme of $\bP(V') \times H'$, finite over $\bP(V)\times H'$
via the projection $\bP(V') - \to \bP(V)$. The corresponding
diagonal embedding
$\cX' \hookrightarrow \bP(V') \times \bP(V) \times H'$
is $A$-equivariant. Thus, it yields an $A$-morphism 
$$H' \to \Hilb(\bP(V')\times \bP(V)).$$
Together with \cite[Theorem 6.3, Proposition 7.2.3]{Ko}, this yields
in turn an $A$-morphism
$$H'^{sn} \to \Chow(\bP(V')\times \bP(V))$$
which sends any point of $H'^{sn}$ to the associated cycle of the
corresponding subscheme of $\bP(V')\times \bP(V)$.
Further, the projection $\bP(V')\times \bP(V) \to \bP(V)$ defines a
morphism $\Chow(\bP(V')\times \bP(V)) \to \Chow \bP(V)$,
by \cite[Theorem 6.8]{Ko}. So we obtain an $A$-invariant morphism 
$H'^{sn} \to \Chow \bP(V)$, that is, a morphism
$$H'^{sn}/A \to \Chow \bP(V).$$
But $H'^{sn}/A$ is semi-normal (e.g., by \cite[Lemma 2.1]{AB-I}), and
hence isomorphic to $H^{sn}$. This yields the desired morphism 
$\gamma$.

To show that $\gamma$ is finite, it suffices to check that the cycle
map has finite fibers, as $M$ is proper. Consider a PSSV $(X,f)$
of type $Q$ over $\bP(V)$. Then $f(X)$ is uniquely determined by the
cycle $\gamma(X)$. By Lemma \ref{finitefibers}, it follows that there
are only finitely many possibilities for $(X,f)$.
\end{proof}

\subsection{Group actions}\label{sec:group-actions}

The group $\GL(V)^G$ acts on the scheme $M_{\Gamma,\bP(V)}$
via the natural action of its quotient $\Aut^G\bP(V)$. Recall the
isomorphisms 
$V \simeq \bigoplus_{\lambda \in F} E(\lambda) \otimes V(\lambda)$ and 
$\GL(V)^G \simeq \prod_{\lambda \in F} \GL(E(\lambda))$.
In particular, if the $G$-module $V$ is multiplicity-free, then
$$
V \simeq \bigoplus_{\lambda\in F} V(\lambda) =: V_F,
$$
and $\GL(V)^G$ is the torus $\bG_m^F$ (product of copies of $\bG_m$
indexed by $F$). We then put
$$
M_{\Gamma,F} := M_{\Gamma,\bP(V)}.
$$ 
We now describe the isotropy group $\Stab_{\bG_m^F}(\xi)$ of a closed
point $\xi \in M_{\Gamma,F}$. Let $(X,f)$ be a representative of
$\xi$. The $G$-module spanned by the image of the associated
morphism $\tf: \tX \to V_F$ depends only on $\xi$; let
$F(\xi)\subseteq F$ be its weight set. Then one easily checks that 
$$
\Stab_{\bG_m^F}(\xi) \simeq \bG_m^{F \setminus F(\xi)}
\times \Aut^G(X,L)/\Aut^G(X,f),
$$
where $L = f^*\cO(1)$. 

Returning to an arbitrary $G$-module $V$ with weight set $F$, we will
establish a bijective correspondence between $\GL(V)^G$-orbits in
$M_{\Gamma,\bP(V)}$ and $\bG_m^F$-orbits in $M_{\Gamma,F}$. For this,
choose lines $\ell(\lambda)$ in $E(\lambda)$ for all 
$\lambda \in F$. This yields an injective $G$-module map  
$V_F \hookrightarrow V$ and an injective homomorphism
$\bG_m^F \hookrightarrow \GL(V)^G$, 
that we regard both as inclusions. We also regard $M_{\Gamma,F}$ as a
closed subscheme of $M_{\Gamma,\bP(V)}$; this subscheme is invariant
under the subgroup $P_F$ of $\GL(V)^G$ that stabilizes all the lines
$\ell(\lambda)$. Note that $P_F$ is a parabolic subgroup of
$\GL(V)^G$ containing $\bG_m^F$ as a direct factor, and  
$$
\GL(V)^G/P_F \simeq \prod_{\lambda \in F} \bP(E(\lambda)).
$$
Consider again a closed point $\xi \in M_{\Gamma,\bP(V)}$ with
representative $(X,f)$, and the associated morphism $\tf : \tX \to V$. 
Since $\tX$ is multiplicity-free and the image $\tf(\tX)$ only depends
on $\xi$, the span of this image decomposes uniquely as
$$
\bigoplus_{\lambda \in F(\xi)} 
\ell(\lambda,\xi) \otimes V(\lambda),
$$
where every $\ell(\lambda,\xi)$ is a line in $E(\lambda)$, and
$F(\xi)$ is a subset of $F$. In fact, $F(\xi)$ only depends on the
$\GL(V)^G$-orbit of $\xi$. This implies readily the following
statement:

\begin{lemma}\label{orbits}
Let $\Omega$ be a $\GL(V)^G$-orbit in $M_{\Gamma,\bP(V)}$. Then
$\Omega$ meets $M_{\Gamma,F}$ along a unique $\bG_m^F$-orbit,
$\Omega_F$. Moreover, $\Omega$ is a homogeneous bundle with fiber
$\Omega_F$ over $\prod_{\lambda \in F(\xi)} \bP(E(\lambda))$, where 
$\xi$ is any point of $\Omega_F$.
\end{lemma}

By Theorem \ref{finiteness}, $M_{\Gamma,\bP(V)}$ contains only
finitely many orbits of spherical varieties over $\bP(V)$. Let 
$M_{\Gamma,\bP(V)}^{\main}$ (the \emph{main part} of
$M_{\Gamma,\bP(V)}$) denote the union of the closures of these 
orbits in $M_{\Gamma,\bP(V)}$ and define $M_{\Gamma,F}^{\main}$
similarly. Then
$M_{\Gamma,\bP(V)}^{\main} = \GL(V)^G M_{\Gamma,F}^{\main}$, 
by Lemma \ref{orbits}. Since any orbit closure of a torus contains 
only finitely many orbits, this yields:

\begin{theorem}\label{main}
The main part $M_{\Gamma,\bP(V)}^{\main}$ contains only finitely many
orbits of $\Aut^G \bP(V)$.
\end{theorem}

Moreover, by Theorem \ref{finiteness}, only finitely many subgroups of
$\tLambda$ arise as weight groups of spherical varieties over
$\bP(V)$. It follows that there are only finitely many stable limits
of such spherical varieties, up to isomorphism and action of $\Aut^G \bP(V)$.  

In contrast, the full moduli space $M_{\Gamma,\bP(V)}$ may contain
infinitely many orbits of $\Aut^G \bP(V)$, as shown by Example 2.4(1). 
However, the number of \emph{closed} orbits is always finite. Indeed,
by Lemma \ref{orbits}, it suffices to check that $M_{\Gamma,F}$
contains only finitely many fixed points of $\bG_m^F$. Consider such 
a fixed point $\xi$ with representative $f:X \to \bP(V_F)$. 
Each irreducible component of $f(X)$ is a multiplicity-free subvariety
of $\bP(V_F)$, invariant under $\bG_m^F$. 
By \cite[Corollary 3.3]{AB-III}, there are only finitely many such
varieties. So the same holds for $(X,f)$, by Lemma
\ref{finitefibers}. 

Rather than considering the orbit structure of $M_{\Gamma,\bP(V)}$, we
will introduce a coarser stratification, still invariant under 
$\Aut^G \bP(V)$. For this, we first classify the SSV's $(X,f)$ over
$\bP(V)$ having the same building blocks: the irreducible
$G$-subvarieties $Y$, their line bundles $L_Y$, and their weight
sets~$F_Y$.

\begin{lemma}\label{maps}
The set of SSV's $f : X \to \bP(V)$ with the given building blocks is
fibered over a product of projective spaces $\bP(E(\lambda))$, into
principal homogeneous spaces under a diagonalizable group. Moreover,
the automorphism group of $f:X \to \bP(V)$ is also diagonalizable.
\end{lemma}

The groups, which we construct explicitly in the proof, are analogs of
the groups $H^i(\hbM^*)$, $i=0,1$ of \cite{Al}. In the case where $G$
is a torus, $\hbM^*$ is a complex of tori and hence its cohomology
groups are diagonalizable groups. For an arbitrary $G$ and
multiplicity-free $V$, the groups $\hbM^i$ are diagonalizable, so that
the $H^i(\hbM^*)$ are diagonalizable as well. In the general case
where $V$ is arbitrary, we get a fibration into diagonalizable groups
over a product of projective spaces, and this bears close resemblance
to the semiabelian case of \cite{Al}. 

\begin{proof}
Let $Y_i$ denote the irreducible components of $X$, 
$Y_{ij} = Y_i \cap Y_j$ the double intersections, $Y_{ijk}$
the triple intersections, etc. For each $(Y,L_Y)$ we have the group
$\Aut^G(Y,L_Y)$, which is diagonalizable by Lemma \ref{diago}, acting on
the set $\hFun(Y,L_Y,F_Y)$ of morphisms of polarized varieties
$(Y,L_Y) \to \big(\bP(V),\cO(1) \big)$ such that the image of 
$V^*\to H^0(Y,L_Y)$ has weight set $F_Y$.

We first consider the case where $V = V_F$ is multiplicity-free. Then
each $\hFun(Y,L_Y,F_Y) = \bG_m^{F_Y}$ is also a group and the action is
described by a homomorphism 
$\phi_{Y,L} : \Aut^G(Y,L) \to \hFun(Y,L_Y,F_Y)$ so that 
$a.f = \phi_i(a) f$.
  
Let $\hbM^*$ denote the cone of the homomorphism
$$
C^*(\phi):C^*(\Aut)\to  C^*(\hFun).
$$  
Explicitly:
$$
\hbM^0 = C^0(\Aut) = \oplus_i \Aut^G(Y_i,L_i),
$$ 
$$
\hbM^1 = C^1(\Aut) \oplus C^0(\hFun) =
  \oplus_{i<j} \Aut^G(Y_{ij},L_{ij}) \oplus_i \hFun(Y_i,L_i,F_i),
$$
$$
\hbM^2 = C^2(\Aut) \oplus C^1(\hFun) =
  \oplus_{i<j<k} \Aut^G(Y_{ijk},L_{ijk}) 
  \oplus_{i<j} \hFun(Y_{ij},L_{ij},F_{ij})
$$ 
and the differentials $d^i:\hbM^i \to \hbM^{i+1}$ are of
the form $(d^i_{\Aut}\times \phi_i^{(-1)^i}, d^{i-1}_{\hFun})$. 

Fix one variety $f:X \to \bP(V)$ with the given building blocks. Then
any other variety $(X',f')$ over $\bP(V)$ with the same blocks 
differs from $(X,f)$ by an element of 
$Z^1(\hbM) = \ker (\hbM^1 \to \hbM^2)$. Indeed, $Z^1(\hbM)$ describes
all other ways to glue the $(Y_i,L_i)$ together compatible on the
triple intersections, and the maps $f_i$ twisted correspondingly. The
group $B^1(\hbM) = \im (\hbM^0 \to \hbM^1)$ describes the effect of
changing each $(Y_i,L_i)$ by an automorphism on the gluing and
maps~$f_i$.

Hence, $X'$ differs from $X$ by an element of $H^1(\hbM)$, and the
automorphism group is $H^0(\hbM)$. Note that all $\hbM^i$ are
diagonalizable groups and so are the cohomology groups.
  
Now, consider the general case. Then the isomorphism classes of PSSV's
over $\bP(V)$ are still classified by the first cohomology set
$H^1(\hbM)$ of a complex of sets $\hbM^*$ properly understood: as the
collection of pairs $(a_{ij}, f_i)$ coinciding on intersections, i.e.,
$(a_{ij})$ is a cocycle and 
  \begin{displaymath}
    a_{ij}. f_i|_{X_{ij}} = f_j|_{X_{ij}},
  \end{displaymath}
modulo the action of collections $(a_i)$ on $(a_{ij},f_i)$. The
automorphism group is still $H^0(\hbM)$.

Any map $f : X \to \bP(V)$ defines a surjective homomorphism
  \begin{displaymath}
V^* = \bigoplus_{\lambda\in F} E(\lambda)^* \otimes V(\lambda)^*
    \onto H^0(X,L)_{F(X)} =
    \bigoplus_{\lambda\in F(X)} V(\lambda)
  \end{displaymath}
which is a point of $\prod_{\lambda\in F(X)} \bP(E(\lambda))$. If
one point appears for some $X$ then so do all others. And for a
fixed point $P$, the first cohomology set is a principal homogeneous
space under $H^1(\hbM^*)$, by the multiplicity-free case. So we are
done. 
\end{proof}

\begin{theorem}\label{rational}
The moduli space $M_{\Gamma,\bP(V)}$ has a natural stratification by
locally closed subsets, invariant under $\Aut^G \bP(V)$. Each subset
is fibered over a product of projective spaces, with fibers being
principal homogeneous spaces over a diagonalizable group.
 
In particular, each irreducible component of $M_{\Gamma,\bP(V)}$ is a
rational variety.
\end{theorem}

\begin{proof}
Let $A$ be the set of possible moment sets $Q$, subdivisions $\cQ$ of
$Q$ into moment polytopes of spherical varieties, and weight sets
$F'$. Then $A$ is a finite set and by Remark \ref{remarks}(2) the
corresponding strata $M_{\alpha}$, $\alpha \in A$, are locally closed
subschemes of $M_{\Gamma,\bP(V)}$.   

Next, consider a finer stratification in which, in addition to the
above data, the isomorphism classes of the irreducible components
$(Y_i,L_i)$ are fixed as well. We claim that these finer strata are
locally closed subsets. 

Indeed, let $f\times \pi:(\cX,\cL) \to \bP(V)\times S$ be a family of
SSV's over $\bP(V)$ with base $S$, and consider the sheaf of
$\cO_S$-algebras  $\cR = \cR(\cX,\cL)$. The identities
$m_{\tilde\lambda,\tilde\mu}=0$ in the definition of $M_{\alpha}$
imply that for each polytope $Q_i \in \cQ$, sending to zero the
$\tlambda$-components outside the cone over $Q_i$ gives a quotient 
algebra $\cR_i$ of $\cR$, and hence a closed subfamily of PSV's 
$(\cY_i,\cL_i)$.

The subsets of $S$ where the fibers $(Y_i,L_i)$ are isomorphic are
locally closed, since they correspond to unions of orbits in
the moduli space of spherical subvarieties of an affine space (cf.~the 
proof of Theorem~\ref{finiteness}). Finally, Lemma~\ref{maps} gives
the stated structure of these locally closed subsets.
\end{proof}

\begin{remark}
The image of the cycle map $\gamma : M^{sn} \to \Chow \bP(V)$ is
contained in the multiplicity-free part $\Chow^{mf}\bP(V)$, consisting
of positive combinations of irreducible multiplicity-free
$G$-subvarieties. This closed subset of $\Chow \bP(V)$ has a natural,
but easier, stratification according to the cycle decomposition and
the isomorphism classes of irreducible components. In homological
language, each stratum corresponds to the global sections of a sheaf
describing the components. The map $\gamma$ is compatible with the 
stratifications. But in general it is neither injective nor
surjective.

In the case where $G$ is a torus, our cycle map is a version of
the Chow morphism from the toric Hilbert scheme to the toric Chow
variety, studied in \cite[Section 5]{HS}.
\end{remark}

Finally, we briefly discuss the moduli space of stable spherical pairs
of type $Q$, i.e., $M_{\Gamma,Q} := M_{\Gamma,\bP(V_Q),Q}$, where 
$V_Q = \bigoplus_{\lambda\in \Gamma \cap Q} V(\lambda)^* \otimes V(\lambda)$.
Then $E(\lambda) = V(\lambda)^*$, so that a natural choice for the
line $\ell(\lambda)$ is the highest weight line in $V(\lambda)^*$. 
The corresponding subspace $M_{\Gamma,F}$ of $M_{\Gamma,Q}$ consists
of the isomorphism classes of those $(X,D)$ such that the divisor $D$
is stable under the action of $U$.

\subsection{Examples}

We reconsider the examples discussed in Section 2.4, from the
viewpoint of moduli spaces.

\medskip

\noindent
1) By Proposition \ref{correspondence}, the moduli space of stable
toric pairs of type $Q$ is $M_{\Gamma,F}$, where 
$F = -\Gamma \cap Q$. This space is described in \cite[Section 2]{Al}. 
Its main part $M_{\Gamma,F}^{\main}$ is irreducible. The normalization
of its image under the cycle map $\gamma$ of Section 4.2 is the toric
variety associated with the secondary polytope of the pair 
$(Q,\Gamma \cap Q)$, see \cite[Chapter 8]{GKZ}.

More generally, one easily checks that $M_{\Gamma,\bP(V),Q}^{\main}$ is
irreducible, if $G = T$ is a torus. This space parametrizes the
normalizations of those $T$-orbit closures in $\bP(V)$ having
moment polytope $Q$, and their limits as stable toric varieties over
$\bP(V)$.

\medskip

\noindent
2) The moduli space of stable reductive pairs is described in
\cite[4.5]{AB-II}. The main tool is a bijective correspondence 
between stable reductive varieties and stable toric varieties with a
compatible action of $W$, see [loc.cit., Theorem 2.8]. For stable
reductive varieties, one checks that the main part of each
$M_{\Gamma,\bP(V),Q}$ is irreducible. 

\medskip

\noindent
3) Likewise, if $G = \SL(2)$ then each $M_{\Gamma,\bP(V),Q}^{\main}$
turns out to be irreducible.

\medskip

\noindent
4) Let $G$, $\Gamma$, $V$, $Q$ and $(X,f)$ be as in Example 2.4(4) and
let $\xi$ be the corresponding point of the space $M_{\Gamma,\bP(V),Q}$. 
Then the connected component of $\xi$ in this space does not meet the
main part (parametrizing the spherical varieties and their stable
limits). This yields an example of a connected component of our moduli
space that contains no irreducible variety. By the results of the
present section, this still holds if $V = V(2,0) \oplus V(4,2)$ is
replaced with a direct sum of any number of copies of $V(2,0)$,
$V(4,2)$.

\medskip

\noindent
5) In the case where $V = V(\lambda)$ is a simple $G$-module, our
moduli space $M_{\Gamma,\bP(V)}$ consists of finitely many points, by
Example 2.4(5). Thus, if $V = E(\lambda)\otimes V(\lambda)$ is an
isotypical $G$-module, then $M_{\Gamma,\bP(V)}$ consists of finitely
many $\GL(E(\lambda))$-orbits, all of them being closed.

\section{Generalizations}\label{sec:generalizations}

\subsection{Split tori}

In the case when $G=\bG_m^r$ is a split torus, the existence of a projective
moduli space and most of the other constructions of this paper can be extended
to an arbitrary base scheme, including $\Spec\bZ$. We point out the changes
that have to be made:

For the construction of the moduli space $\cM_{\Gamma,\bP(V)}$ we used the
$\wG$-Hilbert scheme $\Hilb^{\wG}$ constructed in \cite{AB-III} over an
algebraically closed field of characteristic zero. This has to be replaced
with the multigraded Hilbert scheme, which was constructed in \cite{HS} in
full generality over any base. (Similarly to the classical Grothendieck
Hilbert scheme, even the Noetherian condition may be removed but then one
has to work with locally free instead of flat families).

Quite generally, an action of a split torus $\bG_m^r$, resp. a diagonalizable
group $G$, on an affine scheme $\Spec R$ is the same as a grading of $R$ by
$\bZ^r$, resp.  the dual group $\hG$, a finitely generated abelian group, see
f.e.  \cite{SGA3}.  The arguments in the proof of Theorem~\ref{coarse} are
statements about graded algebras, therefore they apply over a general base.

In particular, the scheme $S'$ obtained by extracting the $N$-th root of a
section $s\in \Gamma(S, \cF^N)$ is $\Spec_{\cO_S}\cA$, where $\cA$ is the
graded algebra $\oplus_{i=0}^{N-1}\cF^{-i}$ with the multiplication given by
$s\in \Hom(\cF^{-N},\cO_S) = \Hom(\cO_S, \cF^N)$, and $S$ is the quotient of
$S'$ by the split diagonalizable group $\mu_N$. Similarly, the group $A$
appearing in the proof is a direct product of several copies of $\mu_N$.

\subsection{Non-algebraically closed fields}

We chose to work over an algebraically closed field of characteristic zero, 
to ease understanding and to avoid introducing cumbersome notation.
However, all the results and arguments may be adapted readily to an
arbitrary field $k$ of characteristic zero, and to a \emph{split} reductive
group $G$. 

Moreover, the main statements \ref{coarse}, \ref{DVR} about existence of
projective coarse space hold for a nonsplit group $G$ as well. A family
$\cX \to S$ of, say affine, $G$-varieties is by definition
multiplicity-free for our purposes, if all the geometric fibers 
$X_{\bar s}$ are multiplicity-free.  For constructing the moduli space, one
can work on an \'etale cover $\coprod S_i$ of $S$, and there exist
such covers with $G\times_S S_i$ split over $S_i$. Similarly, the main theorem
is easily generalized to the case of non-split tori over $\bZ$.

In the non-split setting, we note that the statements about stratifications of
the moduli space have to significantly modified, as some strata glue together.
Also, rationality Theorem~\ref{rational} no longer holds.

\subsection{Stable spherical varieties over a closed subscheme of $\bP(V)$} 

We still fix a subgroup $\Gamma $ of $\tLambda$, and a finite-dimensional
$G$-module $V$. Let $Z$ be a $G$-invariant closed subscheme of $\bP(V)$. 
(We could consider an arbitrary closed subscheme $Z$ but this would not add
any flexibility, since it can be replaced by its unique maximal $G$-invariant
closed subscheme.) 

\begin{lemma}\label{subfunctor-overZ}
The subfunctor $\cM_{\Gamma,Z}$ of stable spherical varieties 
$f:X\to \bP(V)$ with weights in $\Gamma$ that factor through $Z$ is a
closed subfunctor of $\cM_{\Gamma,\bP(V)}$.
\end{lemma}

\begin{proof}
Let $I(Z) \subset \Sym(V^*)$ be the homogeneous ideal of $Z$. Then $Z$ is
the zero subscheme of the homogeneous component $I(Z)_n$ for some $n$ (for
example, for any $n \gg 0$). 
Choose such an $n$ and consider a family 
$f\times \pi :\cX \to \bP(V) \times S$ ; put $\cL := f^* \cO(1)$, so that 
$\pi_*(\cL^n)$ is a locally free $\cO_S$-module by Lemma~\ref{pssv}(vi). The
map 
$I(Z)_n \to H^0(\cX,\cL^n)$ yields a morphism of $\cO_S$-modules 
$\varphi : \cO_S \otimes I(Z)_n \to \pi_*(\cL^n)$.
The zero subscheme of $\varphi$ is the locus of $S$ where our family
factors through $Z$. This is a closed subscheme $S_Z \subseteq S$, whose
formation commutes with base change (since this holds for $\pi_*(\cL^n)$,
see \cite[Proposition III.9.3]{Ha}). 
\end{proof}

\begin{corollary}
The functor $\cM_{\Gamma,Z}$ has a coarse moduli space which is a projective
scheme.
\end{corollary}

Next we obtain an interpretation of the main part of the moduli space
$M_{\Gamma,Z}$, in the case of a torus $T$. Then $Z$ is the disjoint union
of strata $Z_S$, where $S$ denotes a finite subset of $\Lambda$, and  $Z_S$
consists of those points $z = [v] \in Z \subseteq \bP(V)$ such that 
the weight set of $v$ is exactly $S$. Each $Z_S$ is locally closed and
$T$-invariant, and the $T$-action on $Z_S$ factors through a free action of
a quotient torus $T_S$. This stratification is a refinement of the one by
thin Schubert-like cells considered in \cite{Hu}. 

We now make a very simple observation. A point of the quotient 
$Z_S/T = Z_S/T_S$ is the same as the image of a $T$-morphism $f:T_S\to Z_S$,
and the same as the closure $\overline{f(T_S)}\subseteq Z$, which is a
possibly non-normal toric variety. Further, $f$ can be extended in 
a unique way to a finite birational morphism $f : X \to Z$ from the polarized
$T_S$-toric variety $(X,L)$ with moment polytope $Q = \Conv S$. Hence, a point
of $Z_S/T_S$ gives a point of the moduli space $M = M_{\Lambda_S,Z,S}$, where
$\Lambda_S$ is the character group of $T_S$, i.e., the subgroup of
$\Lambda$ generated by differences of elements of $S$.  This observation
may be refined as follows:

\begin{lemma}\label{thin}
The open stratum $M^{\irr}$ of $M$ parametrizing irreducible varieties is
isomorphic to $Z_S/T_S$. Thus, the main part of $M$ is a compactification
of $Z_S/T_S$.
\end{lemma}

\begin{proof}
The polarized toric variety associated with any point of $M^{\irr}$ is
$(X,L)$, and the corresponding finite map to $Z$ restricts to an immersion
on the open orbit $X_0 \simeq T_S$. Thus, each automorphism group 
$\Aut^T(X,f)$ is trivial, and $M^{\irr}$ is a fine moduli space: it admits
a universal family $\cX \to M^{\irr}$. Let $\cX_0$ be the union of all
$T_S$-orbits of maximal dimension. Then we have a morphism $\cX_0 \to Z_S$
and therefore a morphism $M^{\irr} = \cX_0/T_S \to Z_S/T_S$. 

We now construct the inverse of this morphism. Consider the family of tori
$Z_S \to Z_S/T_S$: it can be completed to a family of toric varieties in a
unique way as 
$$
\pi: Z_S \times_{T_S} X \to Z_S/T_S,
$$ 
where $Z_S\times_{T_S} X$ denotes the quotient of $Z_S \times X$ by the
diagonal $T_S$-action. We claim that there exists a $T$-morphism 
$$
f: Z_S \times_{T_S} X \to Z
$$ 
such that the product map 
$\pi \times f: Z_S \times_{T_S} X \to Z_S/T_S \times Z$ is finite, and the
restriction of $f$ to  $Z_S \times_{T_S} X_0 \simeq Z_S$ is the inclusion
of $Z_S$ into $Z$.

To check this, consider the invertible sheaf $\cO(1) \otimes L$ on 
$Z_S \times X$. This sheaf is linearized for the diagonal action of $T_S$,
and hence descends to an invertible sheaf $\cL$ on $Z_S \times_{T_S} X$. 
Moreover, for any $\lambda \in S$, we have a space of global sections
$V^*_{-\lambda}$ of $\cO(1)$, and a line of global sections $k_{\lambda}$
of $L$, which are eigenspaces of $T_S$ of opposite weights. Thus, we obtain
a space of global sections $V^*_{-\lambda} \otimes k_{\lambda}$ of
$\cL$. By the definition of $Z_S$, the direct sum of these subspaces 
(over all $\lambda\in S$) is base-point-free and hence yields a morphism
$$
f: Z_S \times_{T_S} X \to 
\bP(\bigoplus_{\lambda\in S} V_{\lambda} \otimes k_{-\lambda})
\hookrightarrow \bP(V).
$$
One easily checks that the map $\pi \times f$ is finite, and 
$f(Z_S \times_{T_S} X_0) \subseteq Z_S$. Thus, the image of $f$ is
contained in $Z$; this completes the proof of the claim.

By this claim, we obtain a family of $T_S$-toric varieties over $Z$ with
base $Z_S/T_S$, and hence a morphism $Z_S/T_S \to M^{\irr}$ which is the
desired inverse.
\end{proof}

\begin{example}\label{ex:Lafforgue}
  As an application, we give another view to the compactified spaces
  constructed by Lafforgue in \cite{Lafforgue}. In this example, one takes
  $Z$ to be the grassmannian $\Gr$ embedded 
  into a projective space $\bP(V)$ by Pl\"ucker coordinates. Special to this
  case is the fact that a $T$-orbit in $\Gr$ corresponds to a very particular
  ``matroid'' polytope $Q$ (called in \cite{Lafforgue} ``entier''), and the
  weight set $S$ is the full set of lattice points in $Q$. (These polytopes
  have a number of nice properties, for example their integral points
  generate the lattices they lie in.) Thus, the compactification of
  $Z_S/T$ lies already in the open part of the moduli space of stable toric
  varieties over $\bP(V)$ corresponding to subdivisions of $(Q,S)$ into
  matroid polytopes.

Now let us explain this in more detail.
Fix $r\in \bN$ and a free graded $\bZ$-module 
$E=E_0 \oplus E_1 \oplus \cdots \oplus E_n$. The grassmanian $\Gr = \Gr^{r,E}$
parameterizing locally free quotients of $E$ of corank $r$ is a smooth
projective scheme over $\bZ$. It is contained in $\bP(V)$, where 
$$
V = \Lambda^r E = \bigoplus_{\ui=(i_0,\dotsc, i_n),\ \sum i_{\alpha}=r} 
\Lambda^{\ui}E_{\bullet}, 
\quad
\Lambda^{\ui}E_{\bullet} = 
\bigotimes_{\alpha} \Lambda^{i_{\alpha}} E_{\alpha}.
$$
Further, $\Gr$ has a stratification by locally closed strata, usually
called thin Schubert cells 
\begin{displaymath}
    \Gr_{\ud} = \Gr^{r,E}_{\ud} = \{ L \subset E \mid \dim(L\cap E_I) = d_I \}
\end{displaymath}
for some sets of integers $\ud=(d_I)$ labeled by subsets $I$ of
$\{0,1,\dotsc,n\}$ satisfying in particular $d_{\emptyset}=0$,
$d_{\{0,1,\dotsc,n\}} = r$ and $d_I+d_J\le d_{I\cup J} +d_{I\cap J}$. 

Let $\bT=\bG_m^{n+1}$, acting on $E$ in the obvious way, so that
the quotient $T=\bT/\diag(\bG_m)$ acts on $\Gr$ and leaves each $\Gr_{\ud}$
invariant. Note here that if $E$ is multiplicity-free, i.e. all
$\rank E_{\alpha} = 1$, then so is $V$.
Let $S^{r,E}$ be the weight set of the action of $\bT$ on $V$, i.e.,
\begin{displaymath}
    S^{r,E} = \{ 
    (i_0,\dotsc,i_n)\in \bZ^{n+1} \mid
    0\le i_{\alpha} \le \rank E_{\alpha}, \ \sum i_{\alpha} =r 
    \}
\end{displaymath}
and $S^{r,E}_{\bR}=\Conv S^{r,E}$ be the corresponding polytope. Further, let
  \begin{displaymath}
    S^{r,E}_{\ud} = \{ 
    (i_0,\dotsc,i_n)\in S^{r,E} \mid
    \sum_{\alpha\in I}i_{\alpha}\ge d_I 
    \},  \quad
    S^{r,E}_{\ud,\bR} = \Conv(S^{r,E}_{\ud}).
  \end{displaymath}
  We put $S=S^{r,E}_{\ud}$ and $Q=S^{r,E}_{\ud,\bR}$ for simplicity.  By
  \cite[Props.1.1,1.5]{Lafforgue}, $S$ is the common weight set of all
  elements of $\Gr_{\ud}\subset \bP(V)$, and $S= Q\cap \bZ^{n+1}$. Thus
  $\Gr_{\ud} = \Gr_S$ with the preceding notation. Let
\begin{displaymath}
    \oGr_{\ud} = \oGr^{r,E}_{\ud} = \Gr^{r,E}_{\ud} / T .
\end{displaymath}
The main aim of \cite{Lafforgue} is the construction of compactifications
$\oOmega^{r,E}_{\ud}$ of the schemes $\oGr^{r,E}_{\ud}$. The case where all
rank $E_{\alpha} = 1$ and with generic $\ud$ was considered by Kapranov
\cite{Kapranov_ChowQuotients} who used the Chow variety. Another paper
\cite{HKT05} on the multiplicity-free case uses the toric Hilbert
scheme and gives additional moduli interpretation for the 
compactification; see also \cite{Al04}. We will compare Lafforgue's
compactifications with those obtained from Lemma \ref{thin}, that we
denote by $M_{\ud,\Gr}$.

\begin{lemma}
  There exists a finite morphism $\oOmega^{r,E}_{\ud} \to
  M_{\ud,\Gr}$ which restricts to the identity on $\oGr^{r,E}_{\ud}$.
\end{lemma}
\begin{proof}
  The construction of $\oOmega^{r,E}_{\ud}$ is as follows. One first
  constructs a projective morphism of toric varieties $\wcA^S \to \cA^S$
  which is flat with geometrically reduced fibers, such that the generic fiber
  is a toric variety for the polytope $Q$: see
  \cite[Prop.4.3(i)]{Lafforgue}. It immediately follows that this is a family
  of stable toric varieties. 
  
  Recall that the secondary polytope $\Sigma(Q,S)$ is a lattice polytope whose
  faces are in bijection with convex (same as coherent, or regular)
  subdivisions of $Q$ with vertices in $S$, see \cite{GKZ}. Let us denote by
  $\cB^S$ the projective toric variety corresponding to $\Sigma(Q,S)$.  One
  observes that $\cA^S$ is just the open subset of $\cB^S$ corresponding to
  some special convex subdivisions of $(Q,S)$, into matroid polytopes
  $(Q_i,S_i)$ such that $S = \cup S_i$.
  
  Next, one constructs:
  \begin{enumerate}
  \item (\cite[Thm.2.4]{Lafforgue}) A certain scheme, call it $\cA^{S,E}$,
    with a morphism
    \begin{math}
      \cA^{S,E} \to \prod_{\ui\in S} \bP(\Lambda^{\ui}E_{\bullet})
    \end{math}
    such that every fiber is isomorphic to $\cA^S$.
  \item (\cite[Sec.4.3]{Lafforgue}) A family $\pi:\wcA^{S,E} \to \cA^{S,E}$
    which is isomorphic to  $\wcA^{S} \to \cA^{S}$, fiber-wise over 
    $\prod_{\ui\in S} \bP(\Lambda^{\ui}E_{\bullet})$, 
    and a morphism $\wcA^{S,E}\to \bP(V)$ giving a finite morphism $\wcA^{S,E}
    \to \cA^{S,E}\times \bP(V)$.
  \item The closed subscheme $\oOmega^{r,E}$ of $\cA^{S,E}$ corresponding to
    the subfamily mapping to the grassmanian $\Gr\subset \bP(V)$, same as in
    our Lemma~\ref{subfunctor-overZ}.
  \end{enumerate}
  
  Hence, we have a family of stable toric varieties over $\bP(V)$
  parameterized by $\cA^{S,E}$ and a family over $\Gr$ parameterized by
  $\oOmega^{r,E}$, obtained from the previous one by a base change. This gives
  classifying morphisms $f_1: \cA^{S,E} \to M_{\ud,\bP(V)}$ and
  $f_2:\oOmega^{r,E}\to M_{\ud,\Gr}$. To prove that $f_2$ is finite it is
  sufficient to show that $f_1$ is finite to an open subscheme of
  $M_{\ud,\bP(V)}$.
  
  Let $M^0$ be the open part of $M_{\ud, \bP(V)}$ parameterizing varieties
  over $\bP(V)$ with the full weight set $S$. By the mentioned fact about
  moment polytopes of $T$-orbits in $\Gr$, the image of $f_2$ lies in
  $M^0$.  In Section~\ref{sec:group-actions} we showed that $M^0$ has a
  fibration over $\prod_{\ui\in S} \bP(\Lambda^{\ui}E_{\bullet})$
  similar to that above for $\cA^{S,E}$, and every geometric fiber is
  isomorphic to the moduli space $M_{\ud,\bP(V_S)}$ for the multiplicity-free
  module $V_S = \bigoplus_{\ui\in S} \bZ$.
  
  Therefore, the morphism $f_1$ is finite to an open subscheme of
  $M_{\ud,\bP(V)}$ if and only if the morphism from the toric variety $\cA^S$
  to an open subscheme of $M_{\ud,\bP(V_S)}$ is finite for the
  multiplicity-free module $V_S$.  (Thus, we reduced the problem to the case
  when $V$ is multiplicity-free without assuming that $E$ is.) We now prove
  the latter finiteness statement.
  
  The ``Chow toric variety'', or ``Chow toric scheme'' $\cC^S$ is
  defined as the inverse limit of the set of   toric varieties
  corresponding to the fibers of the polytope map from the 
  simplex $\sigma_S$ with vertex set $S$ to $Q$, see 
  \cite{KSZ}. The secondary variety $\cB^S$ is the normalization of the main
  irreducible component of $\cC^S$. One also has a finite ``Chow morphism''
  from the multigraded Hilbert scheme to $\cC^S$ (\cite[Sec.5]{HS}); and by
  our construction of $M_{\ud,\bP(V_S)}$ through the multigraded Hilbert
  scheme, a finite morphism $M_{\ud,\bP(V_S)}\to \cC^S$ (cf. also
  \cite[2.11.11]{Al}. Thus, $\cB^S$ is the normalization of the
  main components of $\cC^S$ and of $\cM_{\ud,\bP(V_S)}$, and the
  statement follows.
\end{proof}
  We expect that $\oOmega^{r,E}_{\ud} \to M_{\ud,\Gr}$ is in fact a closed
  embedding. Note finally that in the multiplicity-free case, the space
  $M_{\ud,\bP(V)}$ is the same as the moduli space of stable toric pairs by
  Section 3.3. The morphism $M_{\ud,\Gr} \to  M_{\ud,\bP(V)}$ in this case
  has an interpretation as the toric analog of the extended Torelli map
  $\overline{M}_g\to \overline{A}_g$, see \cite{Al04}.
\end{example}


\end{document}